\documentclass[12pt,a4paper]{article}

\usepackage[mathscr]{eucal}
\usepackage{enumerate}
\usepackage{epsfig,latexsym,amsfonts,amssymb,amsmath,amscd,
graphics,
epic}
\newcommand{\bp}{\boldsymbol{p}}\newcommand{\bd}{\boldsymbol{d}}\newcommand{\bz}{\boldsymbol{z}}\newcommand{\bn}{\boldsymbol{n}}
\newcommand{\N}{\ensuremath{\mathbb{N}}}
\newcommand{\Z}{\ensuremath{\mathbb{Z}}}

\newcommand{\R}{\ensuremath{\mathbb{R}}}
\newcommand{\C}{\ensuremath{\mathbb{C}}}

\newcommand{\PP}{\ensuremath{\mathbb{P}}}
\newcommand{\redu}
{\ensuremath{\mathrm{red}\:}}\newcommand{\conv}
{\ensuremath{\mathrm{conv}}}\newcommand{\cheb}
{\ensuremath{\mathrm{cheb}}}\newcommand{\jet}
{\ensuremath{\mathrm{jet}}}\newcommand{\Sig}{{\Sigma}}\newcommand{\alp}{{\alpha}}\newcommand{\bet}{{\beta}}
\newcommand{\Int}{\ensuremath{\mathrm{Int}}}

\newcommand{\Idim}{\ensuremath{\mathrm{idim}\:}}
\newcommand{\Sing}{\ensuremath{\mathrm{Sing}}}\newcommand{\eps}{\ensuremath{\varepsilon}}
\newcommand{\Pic}{\ensuremath{\mathrm{Pic}}}

\newcommand{\lcm}{\ensuremath{\mathrm{lcm}}}

\newcommand{\pr}{\ensuremath{\mathrm{pr}}}\newcommand{\del}{{\delta}}\newcommand{\Lam}{{\Lambda}}
\newcommand{\gam}{{\gamma}}

\newcommand{\proofend}{\hfill$\Box$\bigskip}

\newtheorem{theorem}{Theorem}[section]
\newtheorem{lemma}{Lemma}[section]
\newtheorem{proposition}[lemma]{Proposition}
\newtheorem{corollary}[lemma]{Corollary}
\newtheorem{definition}[lemma]{Definition}
\newtheorem{remark}[lemma]{Remark}
\newtheorem{example}[lemma]{Example}

\font\ncsc=cmcsc10  \font\ntt=cmtt12
\begin{document}
\title{Counting curves of any genus on $\PP^2_7$}

\author{Mendy Shoval \and Eugenii Shustin}

\date{}

\maketitle

\begin{abstract}
We compute Gromov-Witten invariants of any genus for del Pezzo
surfaces of degree $\ge2$. The genus zero invariants have been
computed a long ago \cite{DFI,GP}, Gromov-Witten invariants of any
genus for del Pezzo surfaces of degree $\ge3$ have been found by
Vakil \cite{Va}. We solve the problem in two steps: (1) we consider
surfaces $\PP^2_{a,1}$, the plane blown up at $a\ge6$ points on a
conic and one more point outside this conic, and, using techniques
of tropical geometry, obtain a Caporaso-Harris type formula counting
curves of any divisor class and genus subject to arbitrary tangency
conditions with respect to the blown up conic, (2) then we express
the Gromov-Witten invariants of $\PP^2_7$ via enumerative invariants
of $\PP^2_{6,1}$, using Vakil's version of Abramovich-Bertram
formula.
\end{abstract}

\tableofcontents

\section{Introduction}

Del Pezzo surfaces (or Fano surfaces) are smooth rational surfaces
which have an ample anticanonical class $-K$. Each of them is
isomorphic either to the plane $\PP^2$, or to the quadric
$(\PP^1)^2$, or to the plane blown up at $1\le k\le 8$ generic
points (below referred to as of type $\PP^2_k$). The complex
structure of a del Pezzo surface $\Sig$ is generic in the space of
almost complex structures, that particularly implies that the
Gromov-Witten invariants $GW_g(\Sig,D)$ are enumerative, i.e. they
count irreducible complex curve of a given genus $g$ and of a given
divisor class $D$ passing through $-KD+g-1$ generic points on a
surface: for the genus zero it was established in \cite{GP,RT}, for
higher genera observed in \cite[Section 4.2]{Va}. The celebrated
Kontsevich formula \cite{KM} computes recursively all genus zero
Gromov-Witten invariants of the plane. One can find there also
similar formulas for the genus zero Gromov-Witten invariants of
other Del Pezzo surfaces. All genus zero Gromov-Witten invariants of
Del Pezzo surfaces have been computed in \cite{DFI,GP}. In 1998 L.
Caporaso and J. Harris \cite{CH} suggested another formula
recursively computing the Gromov-Witten invariants of the projective
plane for any genus. Vakil \cite{Va} extended the Caporaso-Harris
approach further, producing a formula which computes the number of
curves of a given degree and genus in the plane or on a ruled
surface and, in addition, has few fixed multiple points. In
particular, Vakil computes the Gromov-Witten invariants of arbitrary
genus for the plane blown up at $q\le5$ generic points and the
numbers of curves in a given linear system and of a given genus in
the plane blown up at $6$ points on a conic. In the latter case,
Vakil computes the invariants of $\PP^2_6$ using results of Graber
in Gromov-Witten theory \cite{G} and deriving a suitable extension
of the Abramovich-Bertram formula \cite{AB}. Since then there was no
progress, and here we make the next step calculating the
Gromov-Witten invariants of any genus for $\PP^2_7$, leaving open
the only case of $\PP^2_8$.

The present work goes the way opened in \cite{CH} and developed
further in \cite{Va}. We blow up the plane at $a+1\ge1$ points,
specializing $a$ blown up points on a smooth conic and taking the
$(a+1)$-st blown up point outside, and exhibit a recursive formula
counting curves in any divisor class and of any genus of the
obtained surface $\PP^2_{a,1}$. In the case $a=6$, $\PP^2_{6,1}$ is
almost Fano, and using \cite[Theorem 4.2]{Va}, we convert these
enumerative numbers into genuine Gromov-Witten invariants of
$\PP^2_7$. Our recursive formula is very much similar to the
Caporaso-Harris formula, counting plane curves of any degree and
genus \cite{CH} and to Vakil's formula \cite[Theorem 6.8]{Va}. Its
geometric core is a process of specialization of points in the
constraint configuration (consisting of $-DK+g-1$ points) to the
divisor $E\subset\PP^2_{a,1}$, the strict transform of the conic
through $a$ blown up points. The counted curves degenerate either
into irreducible curves of the same genus with certain tangency
conditions with respect to $E$, or they split off $E$ and few other
components. The comparison of possible degenerations and their
deformations (as the constraint point leaves $E$) yields the
recursion. We should like to mention that for $\PP^2_{a,1}$, one
encounters a {\it new phenomenon}: the reducible degenerate curves
can be {\it non-reduced}. This requires an additional deformation
theory argument based on the local tropical geometry, and this is
the main novelty of the present paper.

Our motivation included an application to the real enumerative
geometry: computation of Welschinger invariants of real Del Pezzo
surfaces of degree $\ge2$ \cite{IKS1,IKS2}, which is based on a
conversion of the complex recursive formula into a formula counting
real rational curves and on the real version of \cite[Theorem
4.2]{Va}. Notice also that $\PP^2_{5,1}$ is just the plane blown up
at $6$ generic points, and our recursive formula directly computes
the Gromov-Witten invariants of $\PP^2_6$ (while in \cite{Va} it is
done indirectly).

\textbf{Acknowledgments.} The second author has enjoyed a support
from the Israeli Science Foundation grant no. 448/09 and from the
Hermann-Minkowski-Minerva Center for Geometry at Tel Aviv
University. A substantial part of this work was done during the
second author's stay at the {\it Mathematisches Forschungsinstitut
Oberwolfach} in the framework of the program {\it Research in
Pairs}. Special thanks are due to I. Itenberg and V. Kharlamov for
very valuable discussions, and to R. Vakil for important remarks.

\section{Counting curves on $\PP^2_{a,1}$}

\subsection{General setting}\label{sec-gs}

{\bf Notation.} Denote by $\Z^\infty_+$ the direct sum of countably
many additive semigroups $\qquad$ \mbox{$\Z_+=\{m\in\Z\ |\ m\ge
0\}$}, labeled by the naturals, and denote by $e_k\in\Z^\infty_+$
the generator of the $k$-th summand, $k\in\N$. For
$\alp=(\alp_1,\alp_2,...)\in\Z^\infty_+$, define
$$
\|\alp\|=\sum_{k\ge1}\alp_k,\quad I\alp=\sum_{k\ge1} k\alp_k,\quad
I^\alp=\prod_{k\ge1} k^{\alp_k},\quad\alp!=\prod_{k\ge1}\alp_k!\ ,
$$ and for
$\alp,\alp^{(1)},\ldots,\alp^{(m)} \in \Z_+^\infty$ such that $\alp
\geq \alp^{(1)} + \ldots + \alp^{(m)}$ we define
$$ \binom{\alp}{\alp^{(1)},\ldots,\alp^{(s)}} =
\frac{\alp!}{\alp^{(1)}! \ldots \alp^{(m)}!(\alp -
\alp^{(1)} - ... - \alp^{(m)})!}\ .$$

\medskip\noindent
{\bf Surfaces under consideration.} Let $\Sig$ be a smooth rational
surface which is a blow-up of $\PP^2$, $E\subset\Sig$ a smooth
rational curve such that
\begin{itemize}\item $-K_\Sig$ is positive on each curve different
from $E$, \item $-(K_\Sig+E)$ is nef and effective.
\end{itemize} Using these data and the genus formula, we immediately derive that
$$-(K_\Sig+E)E=2\ ,$$ that, for each smooth curve $\Lam\ne E$, $$\Lam^2\ge-1\ ,$$ and that, for each smooth rational $(-1)$-curve
$\Lam\ne E$, $$-K_\Sig \Lam=1,\quad E\Lam\le1\ .$$ Thus, the
projection of $\Sig$ onto $\PP^2$ takes $E$ to a straight line or to
a conic, and $\pi:\Sig\to\PP^2$ is a blow-up at several distinct
points such that
\begin{itemize}\item (``planar line-model") either some $a\ge 0$ blown-up points lie on a straight
line $L_0$, $E$ is a strict transform of $L_0$, and $b\le4$ points
are blown-up outside $L_0$, no three blown-up points lie on a line
(except for those on $L_0$), no $6$ blown-up points lie on a conic;
\item (``planar conic-model") or some $a\ge0$ blown up points lie on a smooth conic $C_0$,
$E$ is a strict transform of $C_0$, and $b\le1$ points are blown-up
outside $C_0$, no three blown-up points lie on a line.
\end{itemize}

We will assume that $\Sig$ is {\it generic} in the sense that, in
the above models, the blown-up points are chosen generically subject
to the condition that $a$ of them lie on $L_0$ or $C_0$.

Without loss of generality in our considerations, we can blow up
additional points. Notice that that a planar line-model with $b\ge
3$ can be transformed by Cremona into a planar conic-model. Thus, in
the sequel we will refer only to the planar conic-model with $b=1$
and denote it $\PP^2_{a,1}$. Without loss of generality, we suppose
that $a\ge6$, in particular, $E^2\le-2$.

When referring to the planar conic-model $\pi:\PP^2_{a,1}\to\PP^2$,
we denote by $L$ the pull back of a general line, by $E_1$ the
exceptional divisor disjoint from $E$, by $E_2,...,E_{a+1}$ the
exceptional divisors crossing $E$.

Notice that, in our setting, $(K_\Sig+E)^2=0$, the linear system
$|-(K_\Sig+E)|$ is one-dimensional, its generic element is a smooth
rational curve, and it contains precisely two curves (denoted by
$L',L''$) quadratically tangent to $E$.

Introduce also the semigroups $\Pic_+(\Sig,E)$ and $\Pic(\Sig,E)$
generated by classes $D\in\Pic(\Sig)$ of effective irreducible
divisors such that $DE> 0$ or $DE\ge 0$, respectively. Under our
assumptions $E\not\in\Pic(\Sig,E)$.

\medskip\noindent
{\bf Families of curves.} For a given effective divisor class
$D\in\Pic(\Sig)$ and nonnegative integers $g,n$, denote by
$\overline{\mathcal M}_{g,n}(\Sig,D)$ the moduli space of stable
maps $\bn:\hat C\to\Sig$ of $n$-pointed, connected curves $\hat C$
of genus $g$ such that $\bn_*\hat C\in|D|$ \footnote{In what
follows, by $\bn_*\hat C$ we denote the $\bn$-image of $\hat C$,
whose components are taken with respective multiplicities, by
$\bn(\hat C)$ we denote the reduced image.}. This is an algebraic
stack, whose open dense subset is formed by elements with a smooth
curve $\hat C$ of genus $g$, and the other elements have connected
$\hat C$ of arithmetic genus $g$ with at most nodes as singularities
(cf. \cite[Section 2]{Va}).

Let a divisor class $D\in \Pic(\Sig,E)$, an integer $g$, and two
elements $\alp,\bet\in\Z^\infty_+$ satisfy
\begin{equation}
0\le g\le g(\Sig,D)=\frac{D^2+DK_\Sig}{2}+1,\quad I\alp+I\bet=DE\ .
\label{en1}
\end{equation}
Given a sequence $\bp=(p_{ij})_{i\ge1,1\le j\le\alp_i}$ of
$\|\alp\|$ distinct points of $E$, define ${\mathcal V}_\Sig(D, g,
\alp, \bet, \bp)\subset\overline{\mathcal M}_{g,\|\alp\|}(\Sig,D)$
to be the (stacky) closure of the set of elements $\{\bn:\hat
C\to\Sig,\hat\bp\}$, $\hat\bp=(\hat p_{ij})_{i\ge1,1\le j\le\alp_j}$
a sequence of distinct points of $\hat C$, subject to the following
restrictions:
\begin{itemize}
\item $\hat C$ is smooth, $\bn(\hat p_{ij})=p_{ij}$ for all $i\ge1$, $1\le j\le\alp_i$,
\item $\bn^*(\bn_*C\cap E)$ is the following divisor on $\hat C$: \begin{equation}\bn^*(\bn_*C\cap E)=
\sum_{i\ge1,\ 1\le j\le\alp_i}i\cdot\hat p_{ij}+ \sum_{i\ge1,\ 1\le
j\le\bet_i}i\cdot\hat q_{ij}\ ,\label{ediv}\end{equation} where
$\{\hat q_{ij}\}_{i\ge1,1\le j\le\bet_i}$ is a sequence of
$\|\bet\|$ distinct points of $\hat C\setminus\hat\bp$.
\end{itemize} Put
$$R_\Sig(D,g,\bet)=-D(K_\Sig+E)+\|\bet\|+g-1\ .$$ At last, for a
component $V$ of ${\mathcal V}_\Sig(D,g,\alp,\bet,\bp)$, define the
{\it intersection dimension} $\Idim V$ to be the maximal number of a
priori given generic distinct points of $\Sig$ lying in $\bn_*\hat
C$ for an element $\{\bn:\hat C\to\Sig,\hat\bp\}\in V$.

\subsection{Dimension count and generic elements}\label{secdim}
\begin{proposition}\label{pn1} In the notations of section \ref{sec-gs},
let $D\in\Pic(\Sig,E)$, an integer $g\ge0$, and vectors
$\alp,\bet\in\Z^\infty_+$ satisfy (\ref{en1}), and let ${\mathcal
V}_\Sig(D,g,\alp,\bet,\bp)\ne\emptyset$. Then
\begin{enumerate}\item[(1)] $R_\Sig(D,g,\bet)\ge0$, and each component $V$ of
${\mathcal V}_\Sig(D,g,\alp,\bet,\bp)$ satisfies
\begin{equation}\Idim V\le R_\Sig(D,g,\bet)\ ;\label{enov1}\end{equation}
furthermore, the elements $\{\bn:\hat C\to\Sig,\hat\bp\}\in{\mathcal
V}_\Sig(D,g,\alp,\bet,\bp)$ with smooth $\hat C$ and $\#(\bn(\hat
C)\cap E\setminus\bp)\le n\le\|\bet\|$ form a stratum ${\mathcal
V}_{\Sig,n}(D,g,\alp,\bet,\bp)\subset{\mathcal
V}_\Sig(D,g,\alp,\bet,\bp)$ of intersection dimension
\begin{equation}\Idim{\mathcal V}_{\Sig,n}(D,g,\alp,\bet,\bp)\le-(K_\Sig+E)D+g-1+n\ ,\label{enov11}\end{equation}
\item[(2)] if $D\ne sD_0$ for any $s\ge2$ and any
divisor class $D_0$ such that $-(K_\Sig+E)D_0=0$, and if $V$ is a
component of ${\mathcal V}_\Sig(D,g,\alp,\bet,\bp)$ with $\Idim
V=R_\Sig(D,g,\bet)$, then
\begin{enumerate}\item[(2i)] a generic element $\{\bn:\hat
C\to\Sig,\hat\bp\}\in V$ is an immersion, birational onto its image
$C=\bn(\hat C)$ which is nonsingular along $E$; if, in addition
$E^2\ge-3$, then $C$ is nodal;
\item[(2ii)]if $R_\Sig(D,g,\bet)>0$, then the family $V$ has no base
points outside $\bp$, and the generic curve $C=\bn(\hat C)$ crosses
any a priori given curve $C'\subset\Sig$ transversally outside
$\bp$; furthermore, the image $V_\Sig(D,g,\alp,\bet,\bp)\subset|D|$
of the components of ${\mathcal V}_\Sig(D,g,\alp,\bet,\bp)$ of
intersection dimension $\R_\Sig(D,g,\bet)$ by the projection
\begin{equation}\pr:{\mathcal
V}_\Sig(D,g,\alp,\bet,\bp)\to|D|,\quad\{\bn:\hat C\to\Sig,\hat
\bp\}\mapsto\bn_*\hat C\ ,\label{eproj}\end{equation} is smooth at
$C$.\end{enumerate}\end{enumerate}
\end{proposition}

{\bf Proof.} We divide our argument into several parts: in Steps 1-3
we prove statement (1), in Steps 4-9 we prove statement (2). We
shall use the planar conic-model and the respective notations
introduced in section \ref{sec-gs}.

\smallskip
{\it Step 1: Preliminaries.} Without loss of generality, assume that
${\mathcal V}_\Sig(D,g,\alp,\bet,\bp)$ is irreducible.

From now on and till Step 8 below, we suppose that, for a generic
element $\qquad\qquad$ $\{\bn:\hat C\to\Sig,\hat\bp\}\in {\mathcal
V}_\Sig(D,g,\alp,\bet,\bp)$, the map $\bn:\hat C\to C=\bn(\hat C)$
is birational, {\it i.e.} it is the normalization map.

Then inequalities (\ref{enov1}) and (\ref{enov11}) are evident if
either $D$ is a $(-1)$-curve, or $D=dL-d_1E_1-...-d_{a+1}E_{a+1}$,
where $d\le 2$ or $d_i<0$ for some $i$. So, we suppose that
$D=dL-d_1E_1-...-d_{a+1}E_{a+1}$ with $d\ge 3$, $d_1,...,d_{a+1}\ge
0$.

Due to generic position of $\bp$, one has
$$
\Idim {\mathcal V}_\Sig(D,g,\alp,\bet,\bp)=\Idim {\mathcal
V}_\Sig(D,g,0,\alp+\bet,\emptyset)-\|\alp\|\ ,
$$ and, furthermore, if $\{\bn:\hat C\to\Sig,\hat\bp\}$ is a generic element of ${\mathcal
V}_\Sig(D,g,\alp,\bet,\bp)$, then $\{\bn:\hat C\to\Sig\}$ is generic
for ${\mathcal V}(D,g,0,\alp+\bet,\emptyset)$. Hence we can let
$\alp=0$ and $\bp=\emptyset$. To shorten notation, we write (within
the present proof) ${\mathcal V}(D,g,\bet)$ for ${\mathcal
V}_\Sig(D,g,0,\bet,\emptyset)$.

Inequality (\ref{enov1}) and its proof are completely analogous to
\cite[Propositions 2.1 and 2.2]{CH}: namely, the argument is
developed for curves on any algebraic surface, and its applicability
amounts to checking a number of sufficient numerical conditions,
verified in \cite{CH} for the case of the plane. In the following we
do not copy the reasoning of \cite{CH}, but go through all numerical
conditions and verify them in our setting.

\smallskip
{\it Step 2: Proof of (\ref{enov1}).} The computation of $\Idim
{\mathcal V}(D,g,\bet)$ literally goes along \cite[Section 2.3]{CH},
where one has to verify the following.
\begin{enumerate}
\item[(D1)] The conclusion of \cite[Corollary 2.4]{CH} reads
in our situation as $$\Idim {\mathcal V}(D,g,(DE)e_1)\le
-DK_\Sig+g-1\ ,$$ and it holds, since the hypothesis of
\cite[Corollary 2.4]{CH}, equivalent to $DK_\Sig<0$ is true for any
effective divisor class $D\in\Pic(\Sig,E)$.
\item[(D2)] The inequality
$\deg(\bn^*{\mathcal O}_{\Sig}(-\bd))\ge0$ in \cite[Page 363]{CH},
where
$$\bd=\sum_{i\ge1,\ 1\le j\le\alp_i}i\cdot\hat p_{ij}+ \sum_{i\ge1,\
1\le j\le\bet_i}(i-1)\cdot\hat q_{ij}$$ (cf. (\ref{ediv})), reads in
our situation as $\deg(\bn^*{\mathcal O}_{\Sig}(E)(-\bd))\ge0$, and
it holds, since $$ \deg(\bn^*{\mathcal
O}_{\Sig}(E)(-\bd))=DE-\deg\bd\overset{\text{(\ref{en1})}}{=}\|\bet\|\ge0\
.$$
\item[(D3)] The
inequality \begin{equation}\deg(c_1({\mathcal N}_{\hat
C}(-\bd))\otimes\omega^{-1}_{\hat C})> 0\label{enn5}\end{equation}
in \cite[Page 363]{CH} (${\mathcal N}_{\hat C}$ is the normal sheaf
on $\hat C$ and $\omega_{\hat C}$ is the dualizing bundle) reads in
our setting as
$$\displaylines{
-DK_\Sig+ 2g-2 -\deg\bd + 2 - 2g \cr
=-D(K_\Sig+E)+\|\bet\|=(d-d_1)+\|\bet\|
> 0, }
$$
and it holds true since $d-d_1\ge 1$ as $D$ is represented by a
reduced, irreducible curve with $d\ge 3$.\end{enumerate} After that,
as in the end of the proof of \cite[Proposition 2.1]{CH}, we derive
$$\Idim {\mathcal V}(D,g,\bet)\le\deg(c_1({\mathcal
N}_{\hat C}(-\bd))-g+1=R_\Sig(D,g,\bet)\ ,
$$
which completes the proof of (\ref{enov1}).

\smallskip
{\it Step 3: Proof of (\ref{enov11}).} Again we can assume that
${\mathcal V}_{\Sig,n}(D,g,\alp,\bet,\bp)$ is irreducible and that
$\{\bn:\hat C\to\Sig,\hat\bp\}$ is generic in ${\mathcal
V}_{\Sig,n}(D,g,\alp,\bet,\bp)$. Since $\bn$ is supposed to be
birational on its image, the required statement is analogous to the
inequality in \cite[Corollary 2.7]{CH}, and the proof literally
coincides with the last paragraph of the proof of \cite[Corollary
2.7]{CH}.

\smallskip
{\it Step 4: Immersion.} The sufficient conditions for $\bn:\hat
C\to\Sig$ to be an immersion when $\{\bn:\hat
C\to\Sig,\hat\bp\}\in{\mathcal V}(D,g,\bet)$ is generic, are
\begin{enumerate}\item[(I1)] $\deg(c_1({\mathcal N}_{\hat C}(-\bd))\otimes\omega^{-1}_{\hat
C})\ge 2$, equivalent to $d-d_1+\|\bet\|\ge2$, for the immersion
away from $\bn^{-1}(C\cap E)$ (see \cite[First paragraph of the
proof of Proposition 2.2]{CH}); \item[(I2)] $\deg(c_1({\mathcal
N}_{\hat C}(-\bd))\otimes\omega^{-1}_{\hat C})\ge 4$, equivalent to
$d-d_1+\|\bet\|\ge4$ for the immersion at $\bn^{-1}(C\cap E)$ (see
\cite[Second paragraph of the proof of Proposition 2.2]{CH}).
\end{enumerate}

Suppose that a generic $C=\bn(\hat C)$ has a local singular branch
at $z\in\Sig\setminus E$. Then by B\'ezout's bound for the
intersection of $C$ with a line $l\in|L-E_1|$ passing through $z$,
we get $d\ge d_1+2$, a contradiction by (I1).

Suppose that a generic $C=\bn(\hat C)$ has a singular local branch
at $z\in E$. Then $\|\bet\|\ge1$, and by the above B\'ezout's bound
we get $d-d_1+\|\bet\|\ge3$. Thus, by (I2), it only remains to
analyze the case $d-d_1+\|\bet\|=3$, that means: $d_1=d-2$,
$\|\bet\|=1$ (equivalent to $(C\cdot E)_z=2d-d_2-...-d_{a+1}$), and
$z$ is a center of a unique local branch of multiplicity $2$, {\it
i.e.} a singularity $A_{2s}$, $s\ge1$. Since $\del(A_{2s})=s$, by
the genus formula we have
\begin{equation}s\le\frac{D^2+DK_\Sig+2}{2}-g=d-2-\sum_{i=2}^{a+1}\frac{d_i(d_i-1)}{2}-g\
.\label{egenus}\end{equation} We also recall, that after $s$
successive blow-ups, a curve germ of type $A_{2s}$ turns into a
non-singular branch quadratically tangent to the last exceptional
divisor, and the multiplicities of the exceptional divisors equal
$2$.

Suppose that \begin{equation}2d-d_2-...-d_{a+1}=(C\cdot
E)_z=2s+1\label{einter}\end{equation} (the maximal possible
intersection number). We shall show that, in such a case, the family
$V_1$ of curves of genus $g$ in the linear system $|D|$, having
singularity $A_{2s}$ at $z$, has dimension $\le g$. This will be a
contradiction to the generality of $C$, since, allowing the point
$z$ to move along $E$, we would get
$$\Idim{\mathcal V}(D,g,\bet)\le g+1=R_\Sig(D,g,\bet)-1\ .$$
Indeed, consider the blow-up $\pi:\Sig'\to\Sig$ of $z$ and of $s-1$
more infinitely near points of $C$ at $z$ (all the curves of the
family $V_1$ have the same $s$ blown-up infinitely near points).
Then the considered family of curves on $\Sig$ goes to the family of
immersed curves of genus $g$ in the linear system $|C^*|$ (the upper
asterisk denotes the strict transform) quadratically tangent to the
last exceptional divisor at the intersection point $z^*$ of $C^*$
with $E^*$. We have by (\ref{einter})
$$\deg(c_1(\widetilde{\mathcal
N}_{\hat
C}(-2z^*)))-(2g-2)=-K_{\Sig'}C^*-2=3d-d_1-...-d_{a+1}-2s-2=1>0\ ,$$
(here $\widetilde{\mathcal N}_{\hat C}$ is the normal bundle of
$\widetilde\bn:\hat C\to C^*\hookrightarrow\widetilde\Sig'$) which
ensures (cf. \cite[Page 357]{CH}) that the considered family has
dimension
\begin{equation}\le\deg(c_1(\widetilde{\mathcal N}_{\hat C}(-2z^*)))-g+1=g\ ,\label{edim1}\end{equation} and we are
done, provided (\ref{einter}) holds.

Suppose that \begin{equation}2d-d_2-...-d_{a+1}=(C\cdot
E)_z=2k,\quad 1\le k\le s\ .\label{einter1}\end{equation} Assume
that $k<s$. Since $d_1=d-2$, we derive from B\'ezout that
$d_2,...,d_{a+1}\le2$. Denote by $r_1$ and $r_2$ the number of
$d_i$, $2\le i\le a+1$, equal to $1$ and $2$, respectively. Blowing
down the divisors $E_i$ with $d_i=1$ and ignoring the condition to
pass through the obtained $r$ points on (the image of) $E$, we
obtain a germ of the family $V_1$ of curves of genus $g$ having
$r_2$ double points and a singularity $A_{2s}$ on $E$ so that the
intersection number with $E$ at this singular point equals $2k$.
Fixing the position $z$ of $A_{2s}$ and additionally the position of
$s-1$ infinitely near points obtained in $s-1$ successive blow-ups
of this singularity, we obtain a family $V_2$ of curves of dimension
$\dim V_2\ge\dim V_1-1-(s-k)$. Performing $s$ blow-ups of the
singularity $A_{2s}$, we transform the family $V_2$ into the family
$V_3$ of curves on a surface $\Sig'$ which are quadratically tangent
to the last exceptional divisor in a neighborhood of $z^*$ (the
tangency point of $C^*$). In this situation, taking into account the
genus bound (\ref{egenus}) which yields $r_2+s\le d-2$, we have
$$\deg(c_1(\widetilde{\mathcal N}_{\hat
C}(-z^*)))-(2g-2)=-K_{\Sig'}C^*-1=3d-(d-2)-2r_2-2s-1>0\ ,$$ and
hence
$$\dim V_1\le\dim V_2+(1+s-k)=\dim V_3+(1+s-k)$$
$$\quad\;=\deg(c_1(\widetilde{\mathcal N}_{\hat
C}(-z^*)))-g+1+(1+s-k)$$ $$=2d-2r_2-s-k+1+g\ .\;\;\quad\quad\qquad$$
Restoring the condition to pass through $r_1$ additional points on
$E$ and using the fact that these points are in general position, we
derive that the original family of curves has dimension at most
$$2d-2r_2-s-k+1+g-r_1=(2d-2r_2-r_1-2k)-s+k+1+g\overset{\text{(\ref{einter1})}}{=}-s+k+1+g$$ $$\le
g+1=R_\Sig(D,g,\bet)-1\ .$$

\smallskip

{\it Step 5: Nonsingularity along $E$.} From now on we suppose that
$\bn:\hat C\to\Sig$ is an immersion. In what follows, we argue by
contradiction. Namely, assuming that $C$ is singular at a point of
$C\cap E$, we derive that necessarily $\Idim {\mathcal
V}(D,g,\bet)<R_\Sig(D,g,\bet)$.

Suppose that $\bn$ takes $s\ge2$ points of the divisor $\bd$ to the
same point $z\in E$. Fixing the position of $z$ in $E$, we obtain a
subvariety $U\subset {\mathcal V}(D,g,\bet)$ of dimension $\Idim
U\ge \Idim {\mathcal V}(D,g,\bet)-1$. On the other hand, the same
argument as in the proof of \cite[Proposition 2.1]{CH} gives $$\Idim
U\le h^0(\hat C,{\mathcal N}_{\hat C}(-\bd-\bd')),$$ where
$\bd'=\sum_{i,j\ge1}\hat q_{ij}$. Assuming that $c_1({\mathcal
N}_{\hat C}(-\bd-\bd'))\otimes\omega^{-1}_{\hat C}$ is positive on
$\hat C$ and applying \cite[Observation 2.5]{CH}, we get (cf.
\cite[Page 364]{CH}) $$\Idim {\mathcal V}(D,g,\bet)\le\Idim U+1\le
h^0(\hat C,{\mathcal N}_{\hat C}(-\bd-\bd'))+1$$
$$=\deg(c_1({\mathcal
N}_{\hat C}(-\bd-\bd'))-g+2=-DK_\Sig+2g-2-\deg(\bd+\bd')-g+2$$
$$=-D(K_\Sig+E)+g+\|\bet\|-s<R_\Sig(D,g,\bet)\ .$$ The
positivity required above is equivalent to
\begin{equation}-DK_\Sig-\deg(\bd+\bd')=d-d_1
+\|\bet\|-s > 0\ ,\label{ejan6_4}\end{equation} that holds due to
$\|\bet\|\ge s$, $d>d_1$.

\smallskip

{\it Step 6: Nodality for $E^2\ge-3$ (equivalently $a\le7$).}
Suppose that some distinct points $w_1,w_2,w_3$ of $\hat C$ are
mapped to the same point $z\in\Sig\setminus E$. Fixing the position
of this point, we obtain a subvariety $V\subset {\mathcal
V}(D,g,\bet)$ of dimension $\Idim V\ge\Idim {\mathcal
V}(D,g,\bet)-2$. Thus, provided that $c_1({\mathcal N}_{\hat
C}(-\bd-\bd'))\otimes\omega^{-1}_{\hat C}$ is positive on $\hat C$,
we get $$\Idim{\mathcal V}(D,g,\bet)\le\Idim
V+2\le\deg(c_1({\mathcal N}_{\hat C}(-\bd-\bd'))-g+3$$
$$-DK_\Sig+2g-2-\deg(\bd+\bd')-g+2=-D(K_\Sig+E)+g+\|\bet\|-2$$ $$<-D(K_\Sig+E)+g+\|\bet\|-1=R_\Sig(D,g,\bet)\ ,$$ where $\bd'=w_1+w_2+w_3$.
The required positivity amounts to $-DK_\Sig-\deg(\bd+\bd')=d-d_1
+\|\bet\|-3>0$. The only case, when the latter relation fails is
$\bet=0$ (that is $2d=d_2+...+d_{a+1}$) and $d_1=d-3$. In such a
case, the genus bound leads to
$$\frac{(d-1)(d-2)}{2}\ge\frac{(d-3)(d-3-1)}{2}+\sum_{i=2}^{a+1}\frac{d_i(d_i-1)}{2}+3+g$$
$$\Longrightarrow\quad d^2-3d+2\ge
a\frac{2d}{a}\left(\frac{2d}{a}-1\right)+d^2-7a+18$$
$$\Longrightarrow\quad 2d^2-2ad+8a\le0\ ,$$ and the latter never
holds as $a\le7$, since the discriminant is $9a^2-64<0$.

Suppose that $\bn^{-1}(z)=w_1+w_2$, $w_1\ne w_2\in\hat C$ for some
point $z\in\PP^2_{a,1}\setminus E$, and the two local branches of
$C=\bn(\hat C)$ at $z$ intersect with multiplicity $s\ge2$
(singularity $A_{2s-1}$). In suitable coordinates in a neighborhood
of $z$, $C$ is given by an equation $y^2+2yx^s=0$. The tangent space
to the local equisingular stratum in ${\mathcal O}_{\Sig,z}$ is the
ideal $I=\langle y+x^s,x^{s-1}y\rangle$, and it does not contain
$y$. If $C\cap E\ne \emptyset$, we have the inequality
$d-d_1+\|\bet\|>2$, which comes from $\|\bet\|>0$ and the B\'ezout's
bound $d\ge d_1+2$ for the intersection of $C$ with the line
$l\in|L-E_1|$ through $z$, and this inequality can be rewritten as
$\deg({\mathcal N}_{\hat C}(-\bd-w_1-w_2))>2g-2$. Hence $H^1(\hat
C,{\mathcal N}_{\hat C}(-\bd-w_1-w_2))=0$. The latter relation
immediately implies that there is a section of the bundle ${\mathcal
N}_{\hat C}(-\bd)$ which vanishes at $w_1$ and does not vanish at
$w_2$. Moving $\bn:\hat C\to\Sig$ inside ${\mathcal V}(D,g,\bet)$
along this section, we realize a (one-parameter) deformation of the
local equation $y^2+2yx^s=0$ of $C$, tangent to the line spanned by
the element $y\in{\mathcal O}_{\PP^2_{a,1},z}$, and hence breaking
the equisingularity. This contradiction proves that a generic
$C=\bn(\hat C)$ cannot have tangent local branches. The remaining
option is $C\cap E=\emptyset$. Then $\|\bet\|=0$, and either
$d-d_1>2$, which as above allows one to show that $C$ cannot have
tangent local branches outside $E$, or $d_1=d-2$. In the latter case
we have
\begin{equation}d_2+...+d_{a+1}=2d,\quad d_2,...,d_{a+1}\le2\
.\label{etan}\end{equation} Put $r_j=\#\{i\in[2,a+1]\ :\ d_i=j\}$,
$j=1,2$. Then $r_1+2r_2=2d$. The genus bound
$$\frac{(d-1)(d-2)}{2}\ge\frac{(d-2)(d-3)}{2}+r_2+s$$ yields
$$r_1\ge4+2s\ .$$ Blowing down the divisors $E_i$ for $d_i=1$ and
ignoring the condition to pass through $r_1$ fixed points, we get to
the case $\|\bet'\|=r_1>0$ in which we have derived that the family
of curves with a singularity $A_{2s-1}$ has dimension (here $\Sig'$
is the blown-up $\Sig$)
$$<R_{\Sig'}(D',g,\bet')=d-d_1+g+\|\bet'\|-1=g+r_1+1\ .$$ Restoring the
condition to pass through $r_1$ generic points we obtain a family of
dimension $$<g+1=d-d_1+g-1=R_\Sig(D,g,0)$$ which proves that a
generic element of the original family cannot have singularity
$A_{2s-1}$, $s\ge2$.

\smallskip
{\it Step 7: Absence of base points.} Assuming that
$R_\Sig(D,g,\bet)>0$ and ${\mathcal V}_\Sig(D,g,\alp,\bet,\bp)$ has
a base point $z\in\Sig\setminus\bp$, we obtain $$\Idim {\mathcal
V}_\Sig(D,g,\alp,\bet,\bp)\le h^0(\hat C,{\mathcal N}_{\hat
C}(-\bd-w))\ ,$$ where $\bn(w)=z$. To arrive to a contradiction, we
use the inequality
$$\Idim {\mathcal V}_\Sig(D,g,\alp,\bet,\bp)\le h^0(\hat C,{\mathcal
N}_{\hat C}(-\bd-\bd'))$$
$$=-DK_\Sig+2g-2-\deg(\bd+\bd')-g+1$$
$$=-D(K_\Sig+E)+g+\|\bet\|-2<R_\Sig(D,g,\bet)\ .
$$ This relation is based on the positivity of the bundle $c_1({\mathcal
N}_{\hat C}(-\bd-w))\otimes\omega^{-1}_{\hat C}$ on $\hat C$, which
reduces to the inequality $d-d_1+\|\bet\|>1$. The only case when it
does not hold is $d_1=d-1$, $\|\bet\|=0$, but then $g=0$ and
$R(D,0,0)=d-d_1-1=0$, contrary to the assumptions made. The same
computation with $w\in\hat C$ sent to a tangency point of
$C\setminus\bp$ and $C'$ proves that the intersection
$(C\setminus\bp)\cap C'$ is transversal.

\smallskip
{\it Step 8: Smoothness of the family of curves.} The smoothness of
$V_\Sig(D,g,\alp,\bet,\bp)$ at $C$ is equivalent to
$$H^1(\hat C,{\mathcal N}_{\hat C}(-\bd))=0\quad\Longleftrightarrow\quad d-d_1+\|\beta\|>0\ ,$$
which clearly holds in the considered situation.

\smallskip
{\it Step 9: Multiple covers.} Assume that $\bn:\hat C\to C=\bn(\hat
C)$ is an $s$-fold covering, $s\ge2$, for a generic element
$\{\bn:\hat C\to\Sig,\hat\bp\}\in{\mathcal
V}-\Sig(D,g,\alp,\bet,\bp)$. Then $D=sD_0$, $2-2g=s(2-2g')-r$, where
$g'$ is the geometric genus of $C$, $r$ is the total ramification
multiplicity. Using (\ref{enov11}) for the normalization map
$\bn':C^{\bn}\to C\hookrightarrow\Sig$, we get (cf. \cite[Page
366]{CH} and \cite[Page 62]{Va})
$$\Idim {\mathcal V}_{\Sig,n}(D,g,\alp,\bet,\bp)\le-(K_\Sig+E)D_0+g'-1+\#(C\cap
E\setminus\bp)$$
\begin{equation}\le-\frac{(K_\Sig+E)D}{s}+\frac{g-1-r/2}{s}+n\le-(K_\Sig+E)D+g-1+n\
,\label{enov111}\end{equation} the latter inequality coming from the
nefness of $-(K_\Sig+E)$. Thus, we have proven both (\ref{enov1})
and (\ref{enov11}).

Now, the equality $\Idim {\mathcal
V}_\Sig(D,g,\alp,\bet,\bp)=R(D,g,\bet)$, means the equality in
(\ref{enov111}) with $n=\|\bet\|$, and this leaves the only
possibility $(K_\Sig+E)D_0=0$, equivalent to $d=d_1$, and hence
either $D_0=E_i$, $2\le i\le a+1$, or $D_0=L-E_1$. If $D_0=E_i$,
$2\le i\le a+1$, then $R_\Sig(D,g,\bet)=0$. If $D_0=L-E_1$ and
$R_\Sig(D,g,\bet)>0$, then $C=\bn(\hat C)$ is a generic line in the
one-dimensional linear system $L-E_1$, which, of course, crosses
transversally any a priori given curve $C'$.

The proof of Proposition \ref{pn1} is completed. \proofend

\begin{proposition}\label{pnov2}
Let $D\in\Pic_+(\Sig,E)$, and let $V$ be a component of a nonempty
family ${\mathcal V}_\Sig(D,g,\alp,\bet,\bp)$ such that, for a
generic element $\{\bn:\hat C\to\Sig,\hat\bp\}\in V$, $\bn:\hat C\to
C=\bn(\hat C)$ is an $s$-multiple cover, $s\ge2$. Then
\begin{enumerate}\item[(i)] either $C\subset\Sig$ is a smooth rational
$(-1)$-curve crossing $E$, $\alp=0$, $\bp=\emptyset$, $\bet=e_s$,
$R_\Sig(sC,0,e_s)=0$, and $\bn:\hat C\to C$ has at least two
critical points, one of which is the preimage of $C\cap E$ and has
ramification index $s$;
\item[(ii)] or $D=-s(K_\Sig+E)$, $C\in|-(K_\Sig+E)|$ is either $L'$, or $L''$, $g=0$, $\alp=0$,
$\bp=\emptyset$, $\bet=e_{2s}$, $R_\Sig(D,0,e_{2s})=0$, and
$\bn:\hat C\to C$ has at least two critical points, one of which is
the preimage of $C\cap E$ and has ramification index $s$;
\item[(iii)] or $D=-s(K_\Sig+E)$, $C\in|-(K_\Sig+E)|$ is a line crossing $E$
transversally at two points, $g=0$, $\alp=e_s$, $\bp$ is one of the
points of $C\cap E$, $\bet=e_s$, $R_\Sig(D,0,e_s)=0$, and $\bn:\hat
C\to C$ has two critical points projected onto $C\cap E$;
\item[(iv)] or $D=-s(K_\Sig+E)$, $C\in|-(K_\Sig+E)|$ is a line crossing $E$
transversally at two points, $g=0$, $\alp=0$, $\bp=\emptyset$,
$\bet=2e_s$, $R_\Sig(D,0,2e_s)=1$, and $\bn:\hat C\to C$ has two
critical points of ramification multiplicity $s$ projected onto
$C\cap E$.
\end{enumerate}
\end{proposition}

{\bf Proof.} The statement is straightforward from the consideration
in Step 9 of the proof of Proposition \ref{pn1}. \proofend

\begin{proposition}\label{ptn2}
Let $D\in \Pic_+(\Sig,E)$, $g\ge0$, and $\alp,\bet\in\Z^\infty_+$
satisfy (\ref{en1}), and let $\bp=(p_{ij})_{i\ge1,1\le j\le\alp_i}$
be a sequence of generic distinct points of $E$. Assume that
$R_\Sig(D,g,\bet)=0$, and that ${\mathcal
V}_\Sig(D,g,\alp,\bet,\bp)$ contains an element $\{\bn:\hat
C\to\Sig,\hat\bp\}$ with $\bn$ birational onto its image. Then
\begin{enumerate}\item[(i)] either $D$ is a class of a rational $(-1)$-curve ($E_i$
or $L-E_1-E_i$, $i=2,...,a+1$, in the planar conic-model), $g=0$,
$\bet=e_1$, and $\#{\mathcal V}_\Sig(D,0,0,e_1,\emptyset)=1$;
\item[(ii)] or $D=-(K_\Sig+E)$
and $g=0$, $\bet=e_2$, and $\#{\mathcal V}(D,0,0,e_2,\emptyset)=2$;
\item[(iii)] or $D=-(K_\Sig+E)$, $g=0$,
$\alp=\bet=e_1$, and  $\#{\mathcal V}(D,0,e_1,e_1,\bp)=1$;
\item[(iv)] or $-D(K_\Sig+E)=1$, $g=0$, $I\alp=DE>0$, and
$\#{\mathcal V}_\Sig(D,0,\alp,0,\bp)=1$.
\end{enumerate}
\end{proposition}

{\bf Proof.} Straightforward from the formula
$R(D,g,\bet)=-D(K_\Sig+E)+g+\|\bet\|-1=0$ and the nefness of
$-(K_\Sig+E)$. \proofend

\subsection{Enumerative numbers}\label{sec23}
\begin{definition}\label{d23} Given a divisor class $D\in\Pic(\Sig,E)$, an integer $g\ge0$, vectors
$\alp,\bet\in\Z^\infty_+$ satisfying (\ref{en1}), and a sequence
$\bp=(p_{ij})_{i\ge1,1\le j\le\alp_i}$ of $\|\alp\|$ generic points
of $E$, we define
$$N_\Sig(D,g,\alp,\bet)=\begin{cases}0,\quad &
\text{if}\ V_\Sig(D,g,\alp,\bet,\bp)=\emptyset,\\ \deg
V_\Sig(D,g,\alp,\bet,\bp),\quad & \text{if}\
V_\Sig(D,g,\alp,\bet,\bp)\ne\emptyset, \end{cases}$$
($V_\Sig(D,g,\alp,\bet,\bp)$ introduced in Proposition
\ref{pn1}(2ii)).\end{definition}

This number, of course, does not depend on the choice of $\bp$ and
it counts the curves in $V_\Sig(D,g,\alp,\bet,\bp)$ matching a
generic configuration of $R_\Sig(D,g,\bet)$ points in $\Sig\setminus
E$. In Theorem \ref{t1}, section \ref{sec:formula}, we present a
recursive formula for these numbers. The formula allows one to
compute all these numbers starting with the following initial
values:

\begin{proposition}\label{pini}
Let $D\in \Pic_+(\Sig,E)$, $g\ge0$, and $\alp,\bet\in\Z^\infty_+$
satisfy (\ref{en1}). Then
\begin{enumerate}
\item[(1)] $N(\Sig,-s(K_\Sig+E),0,0,2e_s)=1$ for all $s\ge 1$;
\item[(2)] If $R_\Sig(D,g,\bet)=0$, then $N_\Sig(D,g,\alp,\bet)=0$
except for the following cases:
\begin{enumerate}\item[(i)] $N_\Sig(sD_0,0,0,e_s)=1$ for any divisor class $D$ of
a smooth $(-1)$-curve crossing $E$ and for any $s\ge1$;
\item[(ii)] $N_\Sig(-s(K_\Sig+E),0,0,e_{2s})=2$ for all $s\ge1$;
\item[(iii)] $N_\Sig(-s(K_\Sig+E),0,e_s,e_s)=1$ for all $s\ge1$;
\item[(iv)] $N_\Sig(D,0,\alp,0)=1$ for all divisor classes $D$ with
$-D(K_\Sig+E)=1$, $I\alp=DE$.
\end{enumerate}\end{enumerate}
\end{proposition}

{\bf Proof.} Straightforward from Propositions \ref{pnov2} and
\ref{ptn2}. \proofend

\subsection{Degeneration}\label{sec1}

Let $D\in\Pic(\Sig,E)$, $g\ge0$, and $\alp,\bet\in\Z^\infty_+$
satisfy relations (\ref{en1}) and the inequality
$R_\Sig(D,g,\bet)>0$. Let $\bp=(p_{ij})_{i\ge1,1\le j\le\alp_i}$ be
a sequence of $\|\alp\|$ generic points of $E$, and let $p$ be a
generic point in $E\setminus\bp$. Assume that ${\mathcal
V}_\Sig(D,g,\alp,\bet,\bp)\ne\emptyset$ and introduce
$${\mathcal V}^p_\Sig(D,g,\alp,\bet,\bp)=\big\{\{\bn:\hat
C\to\Sig,\hat\bp\}\in{\mathcal V}(D,g,\alp,\bet,\bp)\ :\
p\in\bn(\hat C)\big\}\ .$$

\begin{proposition}\label{pdeg1}
Let $V$ be a component of ${\mathcal V}_\Sig(D,g,\alp,\bet,\bp)$ of
dimension $R_\Sig(D,g,\bet)$, and let $W$ be a component of
$V\cap{\mathcal V}^p_\Sig(D,g,\alp,\bet,\bp)$ of dimension
$R_\Sig(D,g,\bet)-1$. Then a generic element $\{\bn:\hat
C\to\Sig,\hat\bp\}$ of $W$ is as follows.

(1) Assume that $\bn(\hat C)$ does not contain $E$. Then $\hat C$ is
smooth, and there is $k\ge 1$ such that $W$ is a component of
${\mathcal
V}_\Sig(D,g,\alp+e_k,\bet-e_k,\bp\cup\{p_{k,\alp_k+1}=p\})$.

(2) Assume that $\bn(\hat C)$ contains $E$. Then
\begin{enumerate}\item[(i)] $\hat C=\hat E\cup X\cup Y\cup Z$, where $\bn:\hat E\to E$
is an isomorphism, $X$ is the union of components mapped by $\bn$ to
curves positively intersecting with $E$ and not belonging to the
linear system $|-(K_\Sig+E)|$, $Y$ is the union of components mapped
by $\bn$ to lines belonging to the linear system $|-(K_\Sig+E)|$,
$Z$ is the union of components contracted by $\bn$ to points;
\item[(ii)] $X=\bigcup_{1\le i\le k}\hat C^{(i)}$, $Y=\bigcup_{k<i\le m}\hat C^{(i)}$,
where $0\le k\le m$, and $\hat C^{(1)},...,\hat C^{(m)}$ are
disjoint smooth curves; \item[(iii)] for each $i=1,...,m$, the map
$\bn:\hat C^{(i)}\to\Sig$ represents a generic element in a
component of some family ${\mathcal
V}_\Sig(D^{(i)},g^{(i)},\alp^{(i)},\bet^{(i)},\bp^{(i)})$ of
intersection dimension $R_\Sig(D^{(i)},g^{(i)},\bet^{(i)})$ and such
that

- $\sum_{i=1}^mD^{(i)}=D-E$,

-
$\sum_{i=1}^mR_\Sig(D^{(i)},g^{(i)},\bet^{(i)})=R_\Sig(D,g,\bet)-1$,

- $\bp^{(i)}$, $i=1,...,m$, are disjoint subsets of $\bp$,

- $\bn:\hat C^{(i)}\to\Sig$ is a birational immersion onto the image
for each $i=1,...,k$; the images $\bn(\hat C^{(i)})$, $i=1,...,m$,
are nonsingular along $E$, intersect each other transversally and
only at their nonsingular points lying outside $E$;

- any quadruple $(D',g',0,\bet')$ with $R_\Sig(D',g',\bet')=0$
appears in the list $$(D^{(i)},g^{(i)},\alp^{(i)},\bet^{(i)}),\quad
1\le i\le k\ ,$$ at most once; \item[(iv)] if, for some
$i=k+1,...,m$, $\bn(\hat C^{(i)})=L'$ (or $L''$), then $\bn:\hat
C^{(i)}\to L'$ (resp. $L''$) is an isomorphism.
\end{enumerate}
\end{proposition}

{\bf Proof of Proposition \ref{pdeg1}(1,2i,2ii,2iii).} We follow the
lines of \cite[Section 3]{CH}, where a similar statement
(\cite[Theorem 1.2]{CH}) is proven for the planar case, and also use
some arguments from the proof of \cite[Theorem 5.1]{Va}. Again we do
not copy all the details, but explain the main steps and perform all
necessary computations, which, in fact, are consequences of the
bounds (\ref{enov1}) and (\ref{enov11}), and of the relation
$(K_\Sig+E)E=-2$. The truly new claim, which we provide with a
complete proof in the very end, is that, in the case (2), $\bn_*\hat
C$ does not contain multiple $(-1)$-curves.

\smallskip

Consider a generic one-parameter family in the component $V$ of
${\mathcal V}_\Sig(D,g,\alp,\bet,\bp)$ with the central fibre
belonging to $W$. Using the construction of \cite[Section 3]{CH},
one can replace the given family with a family $\{\bn_t:\hat
C_t\to\Sig,\hat\bp_t\}_{t\in(\C,0)}$ having the same generic fibres
and a semistable central fibre such that (cf. conditions (b), (c),
(e) in\cite[Section 3.1]{CH}):
\begin{itemize}\item the family is represented by a surface
${\mathcal C}$ with at most isolated singularities and two morphisms
$\pi_\Sig:{\mathcal C}\to\Sig$, $\pi_\C:{\mathcal C}\to(\C,0)$, so
that for each $t\ne 0$, $\pi_\Sig:{\mathcal C}_t\overset{\rm
def}{=}\pi_\C^{-1}(t)\to\Sig$ is isomorphic to $\bn_t:\hat
C_t\to\Sig$, \item the fibre $\hat C_0={\mathcal C}_0$ is a nodal
curve,
\item the family is minimal with respect to the above properties.
\end{itemize} Notice that we have disjoint sections \begin{equation}t\in(\C,0)\setminus\{0\}\mapsto\hat p_{ij,t},\ i\ge1,\ 1\le
j\le\alp_i,\quad t\in(\C,0)\setminus\{0\}\mapsto\hat q_{ij,t},\
i\ge1,\ 1\le j\le\bet_i\ ,\label{esec1}\end{equation} defined by
(\ref{ediv}) for each fibre: \begin{equation}\bn_t^*(E\cap\bn_t(\hat
C_t))=\sum_{i\ge1,\ 1\le j\le\alp_i}i\cdot\hat
p_{ij,t}+\sum_{i\ge1,\ 1\le j\le\bet_i}i\cdot\hat q_{ij,t}\
,\label{esec2}\end{equation}
$$\bn_t(\hat p_{ij,t})=p_{ij}\in\bp,\ i\ge1,\ 1\le j\le\alp_i\
.$$ They close up at $t=0$ into some global sections, for which
$\hat p_{ij,0}$, $i,j\ge1$, remain disjoint.

\smallskip

In the case (1), if $\hat C_0$ splits into components $\hat
C^{(1)},...,\hat C^{(m)}$, $m\ge 1$, mapped to curves and some
components mapped to points, then the maps $\bn_0:\hat
C^{(i)}\to\Sig$ represent elements in some ${\mathcal
V}_\Sig(D^{(i)},g^{(i)},\alp^{(i)},\bet^{(i)},\hat\bp^{(i)})$,
$i=1,...,m$, for which we have
$$\sum_{i=1}D^{(i)}=D,\quad\sum_{i=1}^m\|\bet^{(i)}\|\le\|\bet\|-1,\quad\sum_{i=1}^m(g^{(i)}-1)\le
g-m$$ (since some section $q_{kj,t}$ closes up at $p$, and at least
$m-1$ intersection points are smoothed out when deforming $\hat C_0$
into $\hat C_t$, $t\ne0$). Hence (cf. \cite[Page 66]{Va})
$$\sum_{i=1}^m\Idim{\mathcal
V}_\Sig(D^{(i)},g^{(i)},\alp^{(i)},\bet^{(i)},\hat\bp^{(i)})=\sum_{i=1}^m(-(K_\Sig+E)D^{(i)}+g^{(i)}+\|\bet^{(i)}\|-1)$$
$$\le-(K_\Sig+E)D+d+\|\bet\|-1-m=R_\Sig(D,g,\bet)-m\ ,$$ which
implies $m=1$ and thus, in view of Proposition \ref{pn1}, the
statement (1).

\smallskip

In the case (2), assume that $\hat C_0=\hat E\cup\hat
C^{(1)}\cup...\cup\hat C^{(m)}\cup Z$, where $\bn_0(\hat E)=E$, the
components $\hat C^{(i)}$, $1\le i\le m$, are mapped by $\bn_0$ to
curves, and $\bn_0(Z)$ is finite. We have:
\begin{itemize}\item $\hat E$ splits into components $\hat E_1,...,\hat E_a$ such that
$(\bn_0)_*\hat E_i=s_iE$, $s_i\ge1$, $i=1,...,a$;
\item $\{\bn_0:\hat C^{(i)}\to\Sig,\hat\bp^{(i)}\}\in{\mathcal
V}_\Sig(D^{(i)},g^{(i)},\alp^{(i)},\bet^{(i)},\bp^{(i)})$, where
$\bp^{(i)}=\bn_0(\hat\bp^{(i)})\subset\bp$, and $\hat\bp^{(i)}=(\hat
p_{kj,0}\in\hat C^{(i)})_{k,j\ge1}$; furthermore,
$\bet^{(i)}=\widetilde\gam^{(i)}+\gam^{(i)}$, where $\gam^{(i)}$
labels the multiplicities of the sections $\hat q_{kj,t}$ which land
up on $\hat C^{(i)}$ for $t=0$, $i=1,...,m$; \item denote also by
$g'$ the sum of the geometric genera of all the algebraic components
of $Z$, and by $n'$ the number of those connected components $Z'$ of
$Z$, for which $\#Z'\cap\Sing(\hat C_0)$ is greater than the number
of algebraic components of $Z'$ plus $\#\{\hat z\in
Z'\cap\bigcup_i\hat C^{(i)}\ :\ \bn_0(\hat z)\in E\}$.
\end{itemize} By construction, the local branches $\bn_0:(\hat C^{(i)},\hat p_{kj,0})\to\Sig$
and $\bn_0:(\hat C^{(i)},\hat q_{kj,0})\to\Sig$ are deformed
continuously into respective branches $\bn_t:(\hat C_t,\hat
p_{kj,t})\to\Sig$ and $\bn_t:(\hat C_t,\hat q_{kj,t})\to\Sig$,
$t\ne0$. In turn, the other local branches $\bn_0:(\hat C^{(i)},\hat
z)\to\Sig$ with $\bn_0(\hat z)=z\in E$ are deformed so that they
become disjoint from $E$ as $t\ne0$; hence such a point $\hat z\in
\hat C^{(i)}$ must be either an intersection point of $\hat C^{(i)}$
with $\hat E$, or an intersection point of $\hat C^{(i)}$ with a
connected component of $Z$ which joins $\hat C^{(i)}$ with $\hat E$,
and, in the deformation $(\bn_t:\hat C_t\to\Sig)_{t\in(\C,0)}$, this
intersection point is smoothed out. Moreover, all nodes of $\hat
C_0$ are smoothed out in the deformation to $\hat C_t$, $t\ne0$.
Thus, $$\sum_{i=1}^a(g(\hat
E_i)-1)+\sum_{i=1}^m(g^{(i)}-1)+g'+n+n'+\sum_{i=1}^m\|\widetilde\gam_i\|\le
g-1\ ,$$ where $n$ is the number of nodes of $\hat C^{(1)},...,\hat
C^{(m)}$. Applying (\ref{enov11}), we get $$\Idim
W=R_\Sig(D,g,\bet)-1=-(K_\Sig+E)D+g-1+\|\bet\|-1$$
$$\le\sum_{i=1}^m\left(-(K_\Sig+E)D^{(i)}+
g^{(i)}-1+\#(\bn_0(\hat C^{(i)})\cap E\setminus\bp)\right)$$
$$=-(K_\Sig+E)(D-\sum_{i=1}^as_iE)+\sum_{i=1}^m(g^{(i)}-1)+\sum_{i=1}^m\|\bet^{(i)}\|$$
\begin{equation}\le-(K_\Sig+E)D+g-1-\sum_{i=1}^a(2s_i-1+g(\hat
E_i))+\sum_{i=1}^m\|\gam^{(i)}\|-n-n'-g'\
.\label{enov2}\end{equation} Thus, we immediately obtain that
\begin{itemize}\item $a=1$, $s_1=1$, and $g(\hat E_1)=0$, {\it i.e.}
$\bn_0:\hat E\to E$ is an isomorphism; \item
$\sum_{i=1}^m\widetilde\gam^{(i)}=\bet$, {\it i.e.} all the points
$\hat q_{ij,0}$ belong to $\hat C^{(1)}\cup...\cup\hat C^{(m)}$;
\item for each $i=1,...,m$, $\{\bn_0:\hat C^{(i)}\to\Sig,\hat
\bp^{(i)}\}$ is a generic element of a component of ${\mathcal
V}_\Sig(D^{(i)},g^{(i)},\alp^{(i)},\bet^{(i)},\bp^{(i)})$ of
intersection dimension $R_\Sig(D^{(i)},g^{(i)},\bet^{(i)})$; in
particular, all the sections $t\in(\C,0)\mapsto\hat p_{kj,t}$ and
$t\in(\C,0)\mapsto\hat q_{kj,t}$, $k,j\ge1$, are pairwise disjoint
(cf. condition (d) in \cite[Section 3.1]{CH});
\item $n=0$,
{\it i.e.} $\hat C^{(1)},...,\hat C^{(m)}$ are disjoint.
\item $\|\gam^{(i)}\|>0$ for each $i=1,...,m$; \item $g'=0$, {\it
i.e.} all the algebraic components of $Z$ are rational, and $n'=0$,
{\it i.e.} the connected components of $Z$ are chains of rational
curves, which either join the points $\{\hat z\in \bigcup_i\hat
C^{(i)}\ :\ \bn_0(\hat z)\in E,\ \hat z\ne\hat p_{kj,0},\hat
q_{kj,0},\ k,j\ge1\}$ with $\hat E$ (cf. Assumption (b) in
\cite[Section 3.2]{CH}), or are attached by precisely one point to
$\hat E$ or to $\bigcup_i\hat C^{(i)}$; however, the latter type
connected components should not exist due to the minimality of
${\mathcal C}$.
\end{itemize} In view of Propositions \ref{pn1}, \ref{pnov2}, and \ref{ptn2}, it remains to prove that the points $\hat p_{kj,0}$,
which do not belong to $\hat C^{(1)}\cup...\cup\hat C^{(m)}$, lie on
$\hat E$. Indeed, if some point $\hat p_{kj,0}$ were in
$Z\setminus(\hat E\cup\bigcup_i\hat C^{(i)})$, then some point $\hat
z\in\bigcup_i\hat C^{(i)}\setminus\hat\bp_0$ would have been mapped
to $p_{kj}$, which would have resulted in reducing $1$ in the third
and fourth line of bounds (\ref{enov2}), thus, a contradiction.

\smallskip

The final step is to confirm that $C=(\bn_0)_*\hat C_0$ does not
contain multiple $(-1)$-curves. Let, for instance, $E_{a+1}$ have
multiplicity $s\ge2$ in $C$, and (in the above notations) let $\hat
C^{(k+1)},...,\hat C^{(m)}$ be all the components of $\hat C_0$
mapped onto $E_1$. Since $R_\Sig(D,g,\bet)>0$, we have
$\quad\quad\quad\quad$ $D=dL-d_1E_1-...-d_{a+1}E_{a+1}$, $d\ge 1$,
$d_2,...,d_{a+1}\ge 0$. Thus, $C'=(\bn_0)_*(\hat C^{(1)}\cup...\cup
\hat C^{(k)})$ belongs to the linear system
$|(d-2)L-d_1E_1-(d_2-1)E_2-...-(d_a-1)E_a-(d_{a+1}+s-1)E_{a+1}|$.
That is, $C'$ crosses $E_{a+1}\setminus E$ with multiplicity
$d_{a+1}+s-1$, and as explained above these intersection points
persist in the deformation $(\bn_t:\hat C_t\to\Sig)_{t\in(\C,0)}$.
Hence, $(\bn_t)_*\hat C_t$ must cross $E_{a+1}$ with multiplicity
$\ge d_{a+1}+s-1>d_{a+1}$, thus, a contradiction. \proofend

The proof of Proposition \ref{pdeg1}(2iv) will be given in part
(6)of section \ref{sec253}. In turn the statement of Proposition
\ref{pdeg1}(2iv) will not be used before that.

\begin{lemma}\label{l1}
Pick a set $\overline\bp$ of $n-1$ generic points of $\Sig\setminus
E$. Let $\{\bn:\hat C\to\Sig,\hat \bp\}\in{\mathcal
V}_\Sig(D,g,\alp,\bet,\bp)$ and $C:=\bn_*\hat C$ pass through
$\overline\bp$. Then $\{\bn:\hat C\to\Sig,\hat \bp\}$ satisfy the
conditions of Proposition \ref{pdeg1}(2) and it can be uniquely
restored from $C$. Furthermore, $C$ splits into distinct irreducible
components in the following way\ \footnote{The coefficients $s',s''$
in the formula denote the multiplicities of $L',L''$ in $C$,
respectively.}:
\begin{equation}C=E\cup\bigcup_{i=1}^mC^{(i)}
\cup s'L'\cup s''L''\ ,\label{eC}\end{equation}
\begin{enumerate}
\item[(i)] $C^{(1)},...,C^{(m)}$ do not contain
neither $L'$, nor $L''$, and, furthermore, each $C^{(i)}$ either is
a reduced, irreducible curve, or is $kL(p_{kl})$, where $k\ge 2$,
$p_{kl}\in\bp$, and the (uniquely defined) line
$L(p_{kl})\in|-(K_\Sig+E)|$ contains $p_{kl}$, or is $kL(z)$, where
$k\ge 2$, $z\in\overline\bp$, and the (uniquely defined) line
$L(z)\in|-(K_\Sig+E)|$ contains $z$,
\item[(ii)] the curve $C_{\redu}$ is nonsingular along
$E$ and is immersed; furthermore, it is nodal as
$E^2\ge-3$.\end{enumerate} In addition,
\begin{enumerate}\item[(iii)] for each $i=1,...,m$, $C^{(i)}$ is a
generic element in some
$V_\Sig(D^{(i)},g^{(i)},\alp^{(i)},\bet^{(i)},\bp^{(i)})$,
\item[(iv)] $\sum_{i=1}^mD^{(i)}=D-E+(s'+s'')(K_\Sig+E)$, \item[(v)]
$\sum_{i=1}^mR_\Sig(D^{(i)},g^{(i)},\bet^{(i)})=R_\Sig(D,g,\bet)-1$,
\item[(v)] $\bp^{(i)}$, $i=1,...,m$, are disjoint subsets of $\bp$;
\item[(vi)]
$\overline\bp^{(i)}:=\overline\bp\cap C^{(i)}$, $i=1,...,m$, form a
partition of $\overline\bp$;
\item[(viii)] there is a sequence of vectors $\gam^{(i)}\in\Z_+^\infty$, $i=1,...,m$ such that
\begin{equation}0<\gam^{(i)}\le\bet^{(i)}\quad\text{and}\quad
\bet=\sum_{i=1}^m(\bet^{(i)}-\gam^{(i)})\ .\label{eB}\end{equation}
\end{enumerate}
\end{lemma}

{\bf Proof.} Straightforward from Proposition \ref{pdeg1}. \proofend

\subsection{Deformation}\label{sec5}

\subsubsection{Multiplicities} Having generic elements of ${\mathcal
V}^p_\Sig(D,g,\alp,\bet,\bp)$ as described in Proposition
\ref{pdeg1}, we intend to deform them into generic elements of
${\mathcal V}_\Sig(D,g,\alp,\bet,\bp)$ and compute the following
multiplicities.

Let $D\in\Pic_+(\Sig,E)$, a non-negative integer $g$, and vectors
$\alp,\bet\in \Z_+^\infty$ satisfy (\ref{en1}). Pick a sequence $\bp
= \{p_{ij}\}_{i\geq 1, 1 \leq j \leq \alp_i}$ of $\|\alp\|$ generic
distinct points on $E$. Assume, in addition, that
\begin{equation}n:=R_\Sig(D,g,\bet)>0\quad\text{and}\quad D\ne-s(K_\Sig+E),\ s\ge1\ .\label{eassum}
\end{equation} Then, by Proposition \ref{pnov2}, the components of ${\mathcal V}_\Sig(D,g,\alp,\bet,\bp)$
of intersection dimension $R_\Sig(D,g,\bet)$ have the genuine
dimension $R_\Sig(D,g,\bet)$, and their union birationally projects
by (\ref{eproj}) onto its image $V_\Sig(D,g,\alp,\bet,\bp)$ in
$|D|$. Pick a set $\overline\bp$ of $n-1$ generic points of
$\Sig\setminus E$. Then
\begin{equation}V_\Sig(D,g,\alp,\bet,\bp,\overline\bp)\overset{\text{def}}{=}\left\{C\in
V_\Sig(D,g,\alp,\bet,\bp)\ :\
C\supset\overline\bp\right\}\label{estratum}\end{equation} is
one-dimensional (or empty).

Let $p\in E\setminus\bp$ be a generic point, and let $V$ be a
component of ${\mathcal V}^p_\Sig(D,g,\alp,\bet,\bp)$ of
(intersection) dimension $R_\Sig(D,g,\bet)-1$. Let $\{\bn:\hat
C\to\Sig,\hat\bp\}\in V$ be such that $p\in C=\bn_*\hat C$. Within
the given data, put $V(\overline\bp,C)=\emptyset$ if $C\not\in
V_\Sig(D,g,\alp,\bet,\bp,\overline\bp)$, and put $V(\overline\bp,C)$
to be the germ at $C$ of $V_\Sig(D,g,\alp,\bet,\bp,\overline\bp)$ if
$C\in V_\Sig(D,g,\alp,\bet,\bp,\overline\bp)$. By construction,
$\{\bn:\hat C\to\Sig,\hat\bp\}$ is a generic element of $V$; hence
by Propositions \ref{pnov2} and \ref{pdeg1}, $C$ is nonsingular at
$p$. Take a smooth curve germ $L_p$ transversally crossing $E$ and
$C$ at $p$. Thus, we have a well defined map $\varphi_C:V(\overline
\bp,C)\to L_p$. The multiplicities we are interested in are the
degrees $\deg\varphi_C$ (equal $0$ if $V(\overline
\bp,C)=\emptyset)$.

\begin{proposition}\label{p5}
In the above notations and assumptions of Section \ref{sec5}, the
following holds.

(1) Let $\bet_k>0$ and \mbox{$C\in
V_\Sig(D,g,\alp+e_k,\bet-e_k,\bp\cup\{p=p_{k,\alp_k+1}\})$}. Then
$\deg\varphi_C=k$.

(2) Let $C\in|D|$ be given by (\ref{eC}) and satisfy the conditions
of Lemma \ref{l1}. Then
\begin{equation}\deg\varphi_{C}=(s'+1)(s''+1)\sum_{\{\gam^{(i)}\}_{i=1}^m}\prod_{i\in S}I^{\gam^{(i)}}\prod_{i=1}^m
\binom{\bet^{(i)}}{\gam^{(i)}}\ ,\label{e18}\end{equation} where the
sum runs over all sequences $\gam^{(i)}\in\Z_+^\infty$, $i=1,...,m$,
subject to (\ref{eB}), and $$S=\{i\in[1,m]\ :\ C^{(i)}\ne
kL(p_{kl}),\ k\ge2\}\ .$$
\end{proposition}

\subsubsection{Deformation of an irreducible curve}

Here we prove Proposition \ref{p5}(1).

In view of (\ref{eassum}), $C$ possesses the properties listed in
Proposition \ref{pn1}(2). Furthermore, the germ of
$V_\Sig(D,g,\alp+e_k,\bet-e_k,\bp\cup\{p=p_{k,\alp_k+1}\})$ at $C$
is smooth by Proposition \ref{pn1}(2ii). Excluding $p$ from the set
of fixed points, we deduce that $V_\Sig(D,g,\alp,\bet,\bp)$ contains
$C$ and is smooth at $C$ as well. In particular,
$$h^0(\hat C,{\mathcal N}_{\hat C}(-\bd-p))=n-1\quad\text{and}\quad
h^0(\hat C,{\mathcal N}_{\hat C}(-\bd))=n\ .$$ In view of the
generic choice of $\overline\bp$, we obtain
$$h^0(\hat C,{\mathcal N}_{\hat C}(-\bd-p-\overline\bp))=0\quad\text{and}\quad
h^0(\hat C,{\mathcal N}_{\hat C}(-\bd-\overline\bp))=1\ ,$$ which in
its turn means that the germ $V(\overline\bp,C)$ intersect
transversally at $C$ with the hyperplane $H_p=\{C'\in|D|\ :\ p\in
C'\}$ in the linear system $|D|$. The latter conclusion yields that
$V(\overline\bp,C)$ diffeomorphically projects onto the germ of $E$
at $p$ by sending a curve $C'\in V(\overline\bp,C)$ to its
intersection point with $E$ in a neighborhood of $p$. Thus, in
suitable local coordinates $x,y$ of $\Sig$ in a neighborhood of $p$,
we have $p=(0,0)$, $E=\{y=0\}$, $L_p=\{x=0\}$, and a local
parametrization
\begin{equation}C_t=\left\{ay+b(x+t)^k+\sum_{i+kj\ge k}O(t)\cdot
(x+t)^iy^j=0\right\},\quad a,b\in\C^*,\ t\in(\C,0)\
,\label{e30}\end{equation} of $V(\overline\bp,C)$. Hence, for any
point $p'=(0,\tau)\in L_p$, we obtain precisely $k$ curves $C_t\in
V(\overline\bp,C)$ passing through $p'$ and corresponding to the $k$
values of $t$
\begin{equation}t=\left(-\frac{a}{b}\right)^{1/k}\tau^{1/k}+\
\text{h.o.t.}\ ,\label{e15}\end{equation} which completes the proof.

\subsubsection{Deformation of a reducible curve}\label{sec253}

{\bf (1) Geometry of deformation.}

\begin{lemma}\label{l1x} Let $C$ be as in Proposition \ref{p5}(2), and let $V(\overline\bp,C)\ne\emptyset$.
Let $\{C_t\}_{t\in(\C,0)}$ be a parameterized branch of
$V(\overline\bp,C)$ centered at $C=C_0$, and let $\{\bn_t:\hat
C_t\to\Sig,\hat\bp_t\}$, $t\in(\C,0)$ be its lift to ${\mathcal
V}_\Sig(D,g,\alp,\bet,\bp)$. Then in the deformation $C_0\to C_t$,
$t\ne0$, we have
\begin{enumerate}\item[(i)] a point $z\in\Sing(C_{\redu})\setminus E$ is either a node of some $C^{(i)}$, $1\le i\le m$, and then
it deforms into a moving node of $C_t$, $t\ne 0$, or $z$ is a
transverse intersection point of distinct irreducible components
$C',C''$ of $C_{\redu}$, both smooth at $z$, and then the point $z$
deforms into $k'k''$ nodal points of $C_t$, $t\ne0$, where $k'$ and
$k''$ are the multiplicities of $C'$, $C''$ in $C$, respectively;
\item[(ii)] if $z=p_{kl}\in\bp\setminus\bigcup_{i=1}^mC^{(i)}$,
then the germ of $E$ at $z$ turns into a smooth branch of $C_t$
crossing $E$ at $z$ with multiplicity $k$; \item[(iii)] if
$z=p_{kl}\in\bp\cap C^{(i)}$ for some $i=1,...,m$, then the germ of
$C^{(i)}$ at $z$ deforms into a smooth branch of $C_t$ centered at
$z$ and crossing $E$ with multiplicity $k$, and the germ of $E$ at
$z$ deforms into a smooth branch of $C_t$, $t\ne0$, disjoint from$E$
and crossing the former branch of $C_t$ transversally at $k$ points;
\item[(iv)] if $z=q_{kl}\in C^{(i)}\cap E\setminus\bp$ for some
$i=1,...,m$, and $q_{kl}=\bn(\hat q_{kl,0})$, where $(\hat
q_{kl,t})_{t\in(\C,0)}$ is defined by (\ref{esec1}), (\ref{esec2}),
then, first, in case of $C^{(i)}=sL(z)$, $z\in\overline\bp$, we have
$s=k$, and, second, the germ of $C^{(i)}$ at $z$ deforms into a
smooth branch of $C_t$, $t\ne0$, crossing $E$ at one point in a
neighborhood of $z$ with multiplicity $k$, and the germ of $E$ at
$z$ deforms into a smooth branch of $C_t$, $t\ne0$, disjoint from
$E$ and crossing the former branch of $C_t$ transversally at $k$
points;
\item[(v)] if a point $z\in C^{(i)}\cap E\setminus\bp$ is not
$\bn(\hat q_{kl,0})$ for any section $(\hat q_{kl,t})_{t\in(\C,0)}$
as in (iv), then the germ of $E\cup C^{(i)}$ at $z$ deforms into an
immersed cylinder disjoint from $E$;
\item[(vi)] if $s'>0$ (or $s''>0$) then the union of $s'L'$
(resp., $s''L''$) with the germ of $E$ at the point $z=E\cap L'$
(resp., $z=E\cap L''$) turns into an immersed disc disjoint from $E$
and having $s'$ (resp., $s''$) nodes.
\end{enumerate}
\end{lemma}

{\bf Proof.} All the statements follow from Proposition
\ref{pdeg1}(2). \proofend

\medskip

{\bf (2) Strategy of the proof of Proposition \ref{p5}(2).} Let
$$(L_p,p)\to(\C,0),\quad p'\in L_p\mapsto\tau\in(C,0)\ ,$$ be a parametrization of $L_p$.
Assume that $V(\overline\bp,C)\ne\emptyset$, take one of its
irreducible branches $V=\{C_t\ :\ t\in(\C,0)\}$ with a
parametrization normalized by a relation \begin{equation}p'=L_p\cap
C_t\quad\Longleftrightarrow\quad\tau=t^\mu\label{e22}\end{equation}
for some $\mu\ge1$. Clearly, $\deg\varphi_{C}\big|_{V}=\mu$, and the
curves $C\in V$ passing through $p'$ can be associated with
different parameterizations of $V$ obtained by substitutions
$t\mapsto t\eps$, $\eps^\mu=1$. Next, to a parameterized branch $V$
we assign a collection of bivariate polynomials (called {\it
deformation patterns}) which describe deformations of the curve $C$
in neighborhood of isolated singularities of
$\bigcup_{i=1}^mC^{(i)}\cap E$, in neighborhood of non-reduced
components $C^{(i)}$, and in neighborhood of $L',L''$, when moving
along the branch $V$. We show that the constructed deformation
patterns belong to an explicitly described finite set, and then
prove that, given a collection of arbitrary deformation patters from
the above sets, there exists a unique parameterized branch $V$ of
$V(\overline\bp,C)$ which induces the given deformation patterns.
Thus, the degree $\deg\varphi_{C}$ appears to be the number of
admissible collections of deformation patterns in the construction
presented below.

We have chosen deformation patterns in view of their convenience for
further real applications (see \cite{IKS1}), particularly, they suit
well for the calculation of Welschinger signs of real nodal curves.
Another advantage of deformation patterns is that we will be able to
complete the classification of degenerations in Proposition
\ref{pdeg1}, claim (2iv), by showing that the lines $L'$, $L''$
cannot be multiply covered.

In parts (4) and (5) of the present section, we provide a detailed
construction of deformation patterns for isolated singularities and
for multiple lines $L',L''$, which are the most involved cases. In a
similar way, in parts (7)-(9), we construct deformation patterns for
multiply covered lines $L(p_{kl})$ and $L(z)$, but omit routine
details.

The following statement will be used in parts (3)-(9):

\begin{lemma}\label{ms2}
Let $\Pi:\Sig\to\PP^2$ be a contraction of some disjoint
$(-1)$-curves on $\Sig$ which takes $E$ to a straight line or a
conic, and let $x,y$ be coordinates in some affine plane
$\C^2\subset\PP^2$ such that $\C^2\cap \Pi(E)\ne\emptyset$. Then, in
the above notations, the curves $\Pi(C_t)\cap\C^2$ can be defined by
an analytic family of polynomials
\begin{equation}F_t(x,y)=S(x,y)\widetilde G_t(x,y)+t^\mu G_t(x,y),\quad t\in(\C,0)\
,\label{e31ms}\end{equation} where $\{S(x,y)=0\}=\Pi(E)\cap\C^2$,
$\{S(x,y)G_0(x,y)=0\}=\Pi(C_0)\cap \C^2$, and all exponent of $t$
occurring in $G_t(x,y)$ are strictly less than $\mu$.
\end{lemma}

{\bf Proof.} Since $C_0$ contains $E$ as a component, we always can
write $F_t(x,y)=S(x,y)\widetilde G_t(x,y)+t^k G_t(x,y)$ with some
$k>0$ and $\widetilde G_t(x,y)$ not containing term with $t$ to
power $\ge k$ (such terms can simply be moved to $t^k G_t(x,y)$).
Furthermore, since $p$ is a generic point in $E$, we can assume that
$\widetilde G_0(p)\ne0$ and $G_0(p)\ne0$, which together with
(\ref{e22}) yields $k=\mu$. \proofend

\medskip

{\bf (3) Preliminary choice in the construction of deformation
patterns.} First, we fix a sequence $\gam^{(i)}\in\Z_+^\infty$,
$i=1,...,m$, satisfying (\ref{eB}).

Next, for each $i=1,...,m$, take a subset of $\|\gam^{(i)}\|$ points
of $C^{(i)}\cap E\setminus\bp$ selected in such a way that, for
precisely $\gam^{(i)}_k$ of them, the intersection multiplicity of
$C^{(i)}$ and $E$ at the chosen point equals $k$. We notice that the
number of possible choices is
$\prod_{i=1}^m \left(\begin{matrix}\bet^{(i)}\\
\gam^{(i)}\end{matrix}\right)$, and we denote the chosen points by
$q'_{kj}$, $k\ge1$, $1\le j\le\sum_i\gam^{(i)}_k$, assuming that
$(C^{(i)}\cdot E)(q'_{kj})=k$ as $q'_{kj}\in C^{(i)}$. In
particular, for a component $C^{(i)}=kL(p_{kl})$, the only point of
$L(p_{kl})\cap E\setminus\bp$ is always chosen. Denote by $\bz'$ the
set of all the selected points.

In the next parts (4), (5), and (7)-(9), we suppose that
$V(\overline\bp,C)\ne\emptyset$, and there exists a branch $V$ of
$V(\overline\bp,C)$ parameterized with normalization (\ref{e22}) and
which lifts to a family $\{\bn_t:\hat
C_t\to\Sig,\hat\bp_t\}\in{\mathcal V}_\Sig(D,g,\alp,\bet,\bp)$,
$t\in(\C,0)$, such that $(\bn_0)_*\hat C_0=C$, and the points
$\bn_0(\hat q_{ij,0})$ do not belong to $\bz'$ (here the sections
$\hat q_{ij,t}\in\hat C_t$ are uniquely defined by (\ref{esec1}),
(\ref{esec2})). Then put $\bz=\{\bn_0(\hat q_{ij,0})\ :\ i,j\ge1\}$.

\medskip

{\bf (4) Deformation patterns for isolated singularities.} Pick a
point $q'_{kj}\in\bz'\cap \bigcup_{i=1}^mC^{(i)}$.

Let $\Pi:\Sig\to\PP^2$ be the blow down of the disjoint
$(-1)$-curves
$$E_1,...,E_{a-2},\ L-E_{a-1}-E_a,\ L-E_{a-1}-E_{a+1},\
L-E_a-E_{a+1}\ .$$ In particular, it takes $E$ to the line through
the points $\Pi_E(E_2),...,\Pi_E(E_{a-2})$. Take an affine plane
$\C^2\subset\PP^2$ with the coordinates $u,v$ such that
$\Pi_E(E)\cap\C^2=\{v=0\}$. For the sake of notation, we write
$u(z),v(z)$ for the coordinates $u,v$ of the point $\Pi(z)$ as $z$
is a point or a contracted divisor in $\Sig$. For example, we write
$p=(u(p),0)$ with some $u(p)\ne 0$, and, furthermore, we suppose
that $L_p=\{u=u(p)\}$.

Notice that this identification naturally embeds $H^0(\Sig,{\mathcal
O}_\Sig(D))$ into the space $P_{d_0}$ of polynomials in $u,v$ of
degree $\le d_0=2d-d_{a-1}-d_a-d_{a+1}$ (where
$D=dL-d_1E_1-...-d_{a+1}E_{a+1}$).

The curve $C$ is then given by
\begin{equation}F_C(u,v):=v\left(f_1(u)+v\sum_{k,l\ge 0}c_{kl}u^kv^l\right)=0\ ,\label{e20}\end{equation} where
\begin{equation}f_1(u)=\prod_{i=2}^{a-2}(u-u(E_i))^{d_i-1}\prod_{q'_{ij}\in\bz'}(u-u(q'_{ij}))^i\cdot\prod_{z\in\widetilde
C\cap E\setminus\bz'}(u-u(z))^{(\widetilde C\cdot E)(z)}\
,\label{enn10}\end{equation} and $\widetilde C$ is the union of the
components of $C$ different from $E$.

\smallskip

By Lemma \ref{ms2}, the curve $C_t\in V$ is given by an equation
\begin{equation}at^\mu
(f_2(u)+O(t))+v\left(f_1(u)+O(t)+v\sum_{k,l\ge
0}(c_{kl}+O(t))u^kv^l\right)=0\ ,\label{e24}\end{equation} where
\begin{equation}f_2(u)=\prod_{i=2}^{a-2}(u-u(E_i))^{d_i}\prod_{i\ge
1}\prod_{j=1}^{\alp_i}(u-u(p_{ij}))^i\cdot\prod_{i\ge
1}\prod_{j=1}^{\bet_i}(u-u(q_{ij}))^i\label{e21}\end{equation} (here
$p_{ij}$ runs over the points of $\bp$, and $q_{ij}$ are the
projections of the central points $\hat q_{ij,0}$ of the sections
$\hat q_{ij,t}$ defined by (\ref{esec1})), and
\begin{equation}a=-\frac{f_1(u(p))}{f_2(u(p))}\
.\label{e23}\end{equation}

\smallskip

In the coordinates $x=u-u(q'_{ij})$, $y=v$, where $q'_{ij}\in\bz'$
is the chosen point, the curve $C$ is given by an equation
$$c_{02}^{ij}y^2+c_{i1}^{ij}x^iy+\sum_{k+il>2i}c^{ij}_{kl}x^ky^l=0,\quad c^{ij}_{02},c^{ij}_{i1}\ne 0\ ,$$
and the curve $C_t\in V$ ($t\ne 0$) is given by an equation
$$\Phi(x,y,t):=c^{ij}_{02}y^2+c^{ij}_{i1}x^iy+\sum_{k+il=2i}O(t)\cdot x^ky^l+\sum_{k+il>2i}(c^{ij}_{kl}+O(t))x^ky^l$$ $$+
\sum_{k+il<2i}t^{\bn(k,l)}(c^{ij}_{kl}+O(t))x^ky^l=0\ ,$$ where the
last sum does contain nonzero monomials, and for all of them
$\bn(k,l)>0$ and $c^{ij}_{kl}\ne0$, in particular, by (\ref{e24})
and (\ref{e23}),
\begin{equation}\bn(0,0)=\mu,\quad
c^{ij}_{00}=-\frac{f_1(u(p))f_2(u(q'_{ij}))}{f_2(u(p))}\
.\label{e27}\end{equation} Making an additional coordinate change
$$x\mapsto x-yt^{\bn(i-1,1)}\left(\frac{c^{ij}_{i1}}{c^{ij}_{02}}+O(t)\right)(c^{ij}_{i-1,1}+O(t))\ ,$$ we
can annihilate the monomial $x^{i-1}y$ in $\Phi(x,y,t)$.

Notice that
\begin{equation}\bn(k,0)\ge\bn(0,0)=\mu, \quad k>0\ .\label{e25}\end{equation} Indeed, otherwise,
we would have an intersection point of $C_t$ and $E$ converging to
$q'_{ij}$, a contradiction. Hence there exists
\begin{equation}\rho=\min_{k+il<2i}\frac{\bn(k,l)}{2i-k-il}>0\ .\label{e26}\end{equation} Now we
consider the equation
$\Psi(x,y,t):=t^{-2i\rho}\Phi(xt^\rho,yt^{i\rho},t)=0$. We have
$$\Psi(x,y,t)=c^{ij}_{02}y^2+c^{ij}_{i1}x^iy+\sum_{k+il=2i}O(t)\cdot x^ky^l+\sum_{k+il>2i}(c^{ij}_{kl}+O(t))t^{\rho(k+il-2i)}x^ky^l$$
$$+\sum_{k+il\le2i}t^{\bn(k,l)-\rho(2i-k-il)}(c^{ij}_{kl}+O(t))x^ky^l$$
and the well-defined polynomial
$$\Psi(x,y,0)=c_{02}^{ij}y^2+c_{i1}^{ij}x^iy+\sum_{k+il<2i}b^{ij}_{kl}x^ky^l$$ with the
coefficients $$b^{ij}_{kl}=\begin{cases}c^{ij}_{kl},\quad &
\text{if}\ \bn(k,l)=\rho(2i-k-il),\\ 0,\quad & \text{if}\
\bn(k,l)>\rho(2i-k-il)\ .\end{cases}$$ In view of (\ref{e25}),
$b_{k,0}=0$, $k\ge 1$, and hence
\begin{equation}\Psi(x,y,0)=c^{ij}_{02}y^2+c^{ij}_{i1}x^iy+\sum_{k=0}^{i-2}b^{ij}_{k1}x^ky+b^{ij}_{00}\
,\label{e39}\end{equation} where at least one of the $b^{ij}_{k1}$
or $b^{ij}_{00}$ is nonzero. We claim that $b^{ij}_{00}\ne0$.
Indeed, otherwise, $\Psi(x,y,t)=0$ would split into the line
$E=\{y=0\}$ and a curve, crossing $E$ in at least two points, which
in turn would mean that, in the deformation of $C$ into $C_t$ in a
neighborhood $U_{ij}$ of $q'_{ij}$, the germs of $\widetilde C$ and
$E$ at $q'_{ij}$ would glue up with the appearance of at least two
handles, thus breaking the claim of Lemma \ref{l1x}(v). Hence
$b^{ij}_{00}=c^{ij}_{00}$ given by (\ref{e27}), and as observed
above the equation $\Psi(x,y,0)=0$ defines an immersed affine
rational curve. We call $\Psi(x,y,0)$ the {\it deformation pattern}
for the point $q'_{ij}\in\bz'$ associated with the parameterized
branch $V$ of $V(\overline\bp,C)$. Its Newton polygon is depicted in
Figure \ref{f1}. We summarize the information on deformation
patterns in the following statement:

\begin{figure}
\setlength{\unitlength}{1cm}
\begin{picture}(6,4)(-5,0.5)
\thicklines \put(1,1){\line(0,1){2}}\put(1,1){\line(3,1){3}}
\put(1,3){\line(3,-1){3}}
\thinlines\put(1,1){\vector(1,0){4}}\put(1,1){\vector(0,1){3}}\dashline{0.2}(1,2)(4,2)\dashline{0.2}(4,1)(4,2)
\put(3.9,0.5){$i$}\put(0.7,1.9){$1$}\put(0.7,2.9){$2$}
\end{picture}
\caption{Newton polygon of deformation pattern of an isolated
singularity}\label{f1}
\end{figure}
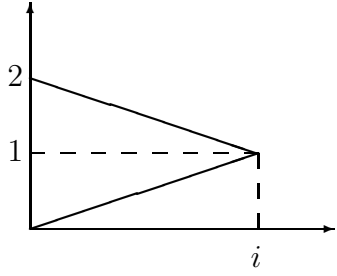

\begin{lemma}\label{l2}
(1) Given a curve $C=E\cup\widetilde C$ and a set $\bz'\subset
\widetilde C\cap E$ as chosen in section \ref{sec253}(3), the
coefficients $c^{ij}_{02}$, $c^{ij}_{i1}$, and
$b^{ij}_{00}=c^{ij}_{00}$ of a deformation pattern of the point
$q'_{ij}\in\bz'$ are determined uniquely up to a common factor.

(2) The set ${\mathcal P}(q'_{ij})$ of polynomials given by formula
(\ref{e39}), having fixed coefficients $c^{ij}_{02}$, $c^{ij}_{i1}$,
and $b^{ij}_{00}=c^{ij}_{00}$  and defining a rational curve,
consists of $i$ elements. Moreover, their coefficients
$b^{ij}_{k1}$, where $k\equiv i+1\mod2$, vanish, and, given a
polynomial $\Psi(x,y)\in{\mathcal P}(q'_{ij})$, the other members of
${\mathcal P}(q'_{ij})$ can be obtained by the following
transformations
\begin{eqnarray}&\Psi(x,y)\mapsto\Psi(x\eps,y),\ \text{where}\
\eps^i=1,\quad\text{if}\ i\equiv1\mod2\ ,\nonumber\\
&\Psi(x,y)\mapsto\Psi(x\eps,y\eps^i),\ \text{where}\
\eps^{2i}=1,\quad\text{if}\ i\equiv0\mod2\ .\nonumber\end{eqnarray}

(3) The germ at $F\in{\mathcal P}(q'_{ij})$ of the family of those
polynomials with Newton triangle $\conv\{(0,0),(0,2),(i,1)\}$, which
define a rational curve, is smooth of codimension $(i-1)$, and
intersect transversally (at one point $\{F\}$) with the affine space
of polynomials
$$\left\{c_{02}^{ij}y^2+c_{00}^{ij}+c_{i1}^{ij}x^iy+\sum_{k=0}^{i-2}c_kx^ky\
:\ c_0,...,c_{i-2}\in\C\right\}\ .$$
\end{lemma}

{\bf Proof.} The first statement comes from construction. In
particular, $a^{ij}_{00}$ is determined by formula (\ref{e27}). The
second statement follows from \cite[Lemma 3.5]{Sh0}. At last, the
second statement additionally implies that the family in the third
statement is the (germ of the) orbit of $F$ by the action
$$F(x,y)\mapsto \xi F(\xi_1x+\xi_0,\xi_2y),\quad\xi,\xi_1,\xi_2\in\C^*,\ \xi_0\in\C\ ,$$ and the third claim follows. \proofend

\begin{remark}\label{r5}
(1) Since the coefficients $b^{ij}_{k1}$ in (\ref{e39}) vanish for
$k\equiv i+1\mod2$, and the other ones not, the values
$$\bn(i-2k,1)=\rho(2i-(i-2k)-i)=\frac{\mu}{2i}\cdot 2k=\frac{\mu}{i},\quad 1\le k\le\frac{i}{2}\
,$$ are integer; hence $\mu$ is divisible by $i$.

(2) The treatment of this section covers the cases of non-multiple
components $C^{(i)}=L(p_{kl})$, $L(z)$, $L'$, or $L''$. In the next
sections we consider nonreduced components $C^{(i)}$, however the
case of multiplicity one is also covered there with the same (up to
a suitable coordinate change) answer.
\end{remark}

\medskip

{\bf (5) Deformation patterns for nonisolated singularities, I.} In
this section we construct deformation patterns for a component
$s'L'$ of $C$ with $s'\ge1$.

\smallskip

{\it Step 1: Preparation.} To relax the notation, within this
section we write $s$ for $s'$.

Consider the blow-down $\Pi:\Sig\to\PP^2$ which contracts
$E_1,...,E_{a+1}$. It takes $E$ to a conic $\Pi(E)$ passing through
$a$ fixed points $\Pi(E_2),...,\Pi(E_{a+1})$, and takes $L',L''$ to
straight lines tangent to $\Pi(E)$ and passing through the fixed
point $\Pi(E_1)$. The linear system $|D|$ on $\Sig$ turns into the
linear system $\Pi_*|D|$ of plane curves of degree $d$ with multiple
points $\Pi(E_1),...,\Pi(E_{a+1})$ of order $d_1,...,d_{a+1}$,
respectively. Take an affine plane $\C^2\subset\PP^2$ with
coordinates $x,y$ such that $\Pi(L')=\{y=0\}$,
$\Pi(E)=\{S(x,y):=y+yx+x^2=0\}$, $\Pi(E_1)$ is the infinite point of
$\Pi(L')$, and the tangency point of $\Pi(L')$ and $\Pi(E)$ is the
origin. By Lemma \ref{ms2}, the curves $\Pi(C_t)$, where $C_t\in V$,
are given by (\ref{e31ms}) in which we can suppose that $\widetilde
G_0(x,y)=y^s\widetilde G'_0(x,y)$ with the polynomial $\widetilde
G'_0(x,y)$ defining the union of the components of $\Pi(C)$
different from $\Pi(L'),\Pi(L'')$ and such that $\widetilde
G'_0(0,0)=1$. Furthermore, the value $c=\widetilde G_0(0,0)\ne0$ can
be computed from relations (\ref{e22}) and (\ref{e31ms}) which yield
\begin{equation}
G_0(\Pi(p))=-\frac{d}{d\tau}\left(S(x,y)
\big|_{\Pi(L_p)}\right)\Big|_{\Pi(p)}\cdot \widetilde G'_0(\Pi(p))\
.\label{e32}\end{equation} 

\smallskip

{\it Step 2: Tropical limit.} Now we find the tropical limit of the
family (\ref{e31ms}) in the sense of \cite{Sh0}. For the reader's
convenience, we shortly recall what is the tropical limit. Consider
the family of the curves $K_t$ defined by polynomials (\ref{e31ms})
for $t\ne0$ in the trivial family of the toric surfaces
$\C(\Delta)$, where $\Delta$ is the Newton polygon of a generic
polynomial in (\ref{e31ms}). Then we close up this family at $t=0$
(all details can be found in \cite{Sh0}):
\begin{itemize}\item to each monomial $x^iy^j$ we assign the point
$(i,j,\nu_{ij})\in\Z^3$, where $\nu_{ij}$ is the minimal exponent of
$t$ in the coefficient of $x^iy^j$ in $F_t$, and then define a
convex piece-wise linear function $\nu:\Delta\to\R$, whose graph is
the lower part of the convex hull of all the points
$(i,j,\nu_{ij})$;
\item the surface $\C(\Delta)$ degenerates into the union of toric
surfaces $\C(\del)$, where $\del$ runs over the maximal linearity
domains of $\nu$ (the subdivision of $\Delta$ into these linearity
domains, which all are convex lattice polygons, we call
$\nu$-subdivision),\;
\item then write
$F_t(x,y)=\sum_{(i,j)\in\Delta}(c_{ij}+O(t))t^{\nu(i,j)}$ and define
the limit of the curves $K_t$ at $t=0$ as the union of the limit
curves $K_\del\subset\C(\del)$ given by limit polynomials
$f_\del=\sum_{(i,j)\in\del}c_{ij}x^iy^j$ for all pieces $\del$ of
the $\nu$-subdivision of $\Delta$.
\end{itemize} Notice that \begin{equation}\nu_{ij}\ge\nu(i,j)\
\text{for all}\ i,j,\quad \text{and}\ \nu_{ij}=\nu(i,j)\ \text{for
the vertices of}\ \nu\text{-subdivision}\ .\label{e37}\end{equation}

\smallskip

\begin{figure}
\setlength{\unitlength}{1cm}
\begin{picture}(13,20)(-2,0)
\thinlines\put(1.5,1){\vector(1,0){4}}\put(1.5,1){\vector(0,1){4.5}}
\put(7,1){\vector(1,0){4}}\put(7,1){\vector(0,1){4.5}}
\put(1.5,7.5){\vector(1,0){4}}\put(1.5,7.5){\vector(0,1){4.5}}
\put(7,7.5){\vector(1,0){4}}\put(7,7.5){\vector(0,1){4.5}}
\put(1.5,14){\vector(1,0){4}}\put(1.5,14){\vector(0,1){4.5}}
\put(7,14){\vector(1,0){4}}\put(7,14){\vector(0,1){4.5}}

\put(4.5,1){\line(0,1){1}}\put(4.5,2){\line(-1,1){3}}
\put(10,1){\line(0,1){1}}\put(10,2){\line(-1,1){3}}
\put(4.5,7.5){\line(0,1){1}}\put(4.5,8.5){\line(-1,1){3}}
\put(10,7.5){\line(0,1){1}}\put(10,8.5){\line(-1,1){3}}
\put(4.5,14){\line(0,1){1}}\put(4.5,15){\line(-1,1){3}}
\put(10,14){\line(0,1){1}}\put(10,15){\line(-1,1){3}}
\put(8,14){\line(0,1){1}}\put(1.5,7.5){\line(1,1){1}}
\put(2,8){\line(1,0){2}}\put(4,7.5){\line(0,1){0.5}}
\put(4,8){\line(1,1){0.5}}\put(7,7.5){\line(2,1){1}}
\put(7,8.5){\line(2,-1){1}}\put(2.5,1){\line(0,1){1}}
\put(2.5,1.5){\line(1,0){2}}\put(7,1.5){\line(2,-1){1}}
\put(8,1){\line(0,1){1}}

\dashline{0.2}(1.5,15)(2.5,15)\dashline{0.2}(2.5,14)(2.5,15)
\dashline{0.2}(8,7.5)(8,8.5)\dashline{0.2}(7,8)(8,8)

\thicklines \put(1.5,2.5){\line(2,-1){1}}\put(2.5,2){\line(1,0){2}}
\put(7,2.5){\line(2,-1){1}}\put(8,2){\line(1,0){2}}
\put(1.5,9){\line(2,-1){1}}\put(2.5,8.5){\line(1,0){2}}
\put(7,9){\line(2,-1){1}}\put(8,8.5){\line(1,0){2}}
\put(1.5,15.5){\line(2,-1){1}}\put(2.5,15){\line(1,0){2}}
\put(7,15.5){\line(2,-1){1}}\put(8,15){\line(1,0){2}}

\put(3.3,0){\rm (e)}\put(8.8,0){\rm (f)} \put(3.3,6.5){\rm
(c)}\put(8.8,6.5){\rm (d)}\put(3.3,13){\rm (a)}\put(8.8,13){\rm (b)}
\put(0.4,15.5){$s+1$}\put(1.1,14.9){$s$}\put(5.9,15.5){$s+1$}
\put(5.9,2.5){$s+1$}\put(6.6,1.4){$1$}\put(5.9,8.4){$j+1$}
\put(6.6,7.9){$j$}\put(2.4,13.5){$2$}\put(7.9,13.5){$2$}
\put(7.9,7){$2$}\put(2.4,0.5){$2$}\put(7.9,0.5){$2$}\put(1.1,17.9){$d$}
\put(6.6,17.9){$d$}\put(1.1,11.4){$d$} \put(6.6,11.4){$d$}
\put(1.1,4.9){$d$} \put(6.6,4.9){$d$}\put(4.2,0.5){$d-d_1$}
\put(9.7,0.5){$d-d_1$}\put(4.2,7){$d-d_1$} \put(9.7,7){$d-d_1$}
\put(4.2,13.5){$d-d_1$} \put(9.7,13.5){$d-d_1$}

\end{picture}
\caption{Deformation patterns for $sL'$}\label{f2}
\end{figure}
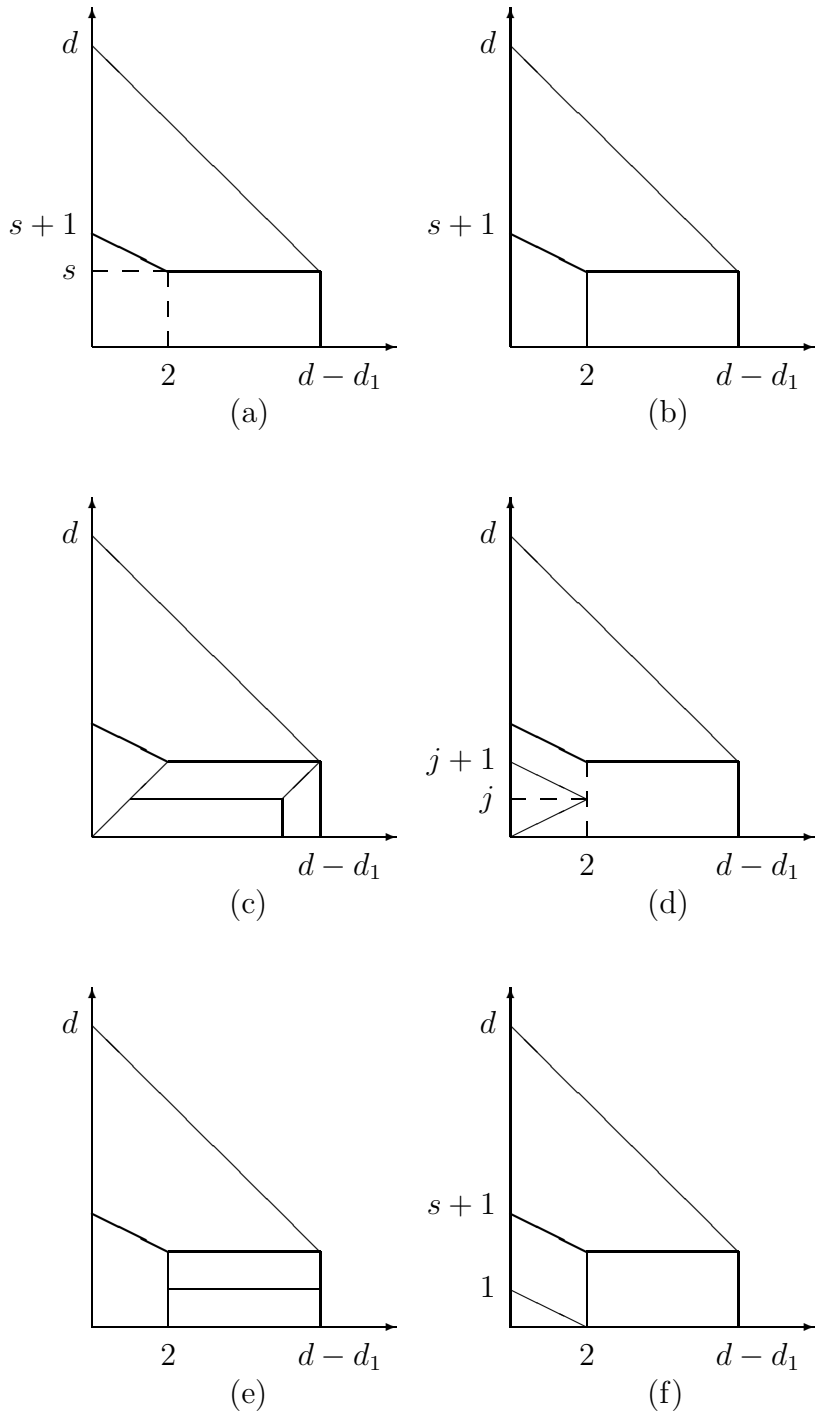

{\it Step 3: Subdivision of the Newton polygon.} In Figure
\ref{f2}(a), we depicted the Newton polygon $\Delta$ of $F_t(x,y)$,
whose part $\del_1$ above the bold line is the Newton polygon of
$S(x,y)y^s\widetilde G'_0(x,y)$, and it is the linearity domain of
$\nu$, where it vanishes. Below the bold line, $\nu$ is positive,
and hence this part $\del_2$ is subdivided by other linearity
domains of $\nu$. We claim that $\del_2$ is subdivided into two
pieces as shown in Figure \ref{f2}(b).

{\it (3A)} The truncation of the polynomial $S(x,y)y^s\widetilde
G'_0(x,y)$ to the incline segment $\sigma_1$ of the bold line is
$(y+x^2)y^s$, and the truncation to the horizontal segment
$\sigma_2$ of the bold line is $x^2y^s\widetilde G'_0(x,0)$. Observe
that the roots of $\widetilde G'_0(x,0)$ are the coordinates of the
intersection points of the closure of $\Pi(C\setminus(L'\cup E))$
with $\Pi(L')$ (counting multiplicities). By Lemma \ref{l1x}(i), the
latter intersection points persist in the deformation $C_t$,
$t\in(\C,0)$, which implies (cf. \cite[Section 3.6]{Sh0}) that the
subdivision of $\del_2$ must contain one or few trapezes $\del$ with
a pair of horizontal sides parallel to $\sigma_2$ of length at least
$\deg\widetilde G'_0(x,0)$, and of the total height $s$ (like in
Figure \ref{f2}(c)), and the corresponding polynomials $f_\del$ are
divisible by $\widetilde G'_0(x,0)$.

{\it (3B)} Let us show that all such trapezes must be rectangles.

First, notice that $\nu(0,0)=\mu$, since
$F_t(0,0)=t^\mu(cG_0(0,0)+O(t))$. The convexity of $\nu$ and its
constancy along the horizontal sides of the trapezes (cf. Figure
\ref{f2}(c)) yield that $0\le\nu\le\mu(s-l)/s$ along such a side on
the height $0\le l<s$. Ordering the monomials of $\widetilde
G_t(x,y)$, first, by the growing power of $x$ and then by the
growing power of $y$, and using (\ref{e31ms}), we inductively derive
that the minimal exponents $\mu_{il}$ in the coefficients of
$x^iy^l$ in $\widetilde G_t(x,t)$ satisfy
\begin{equation}\mu_{il}\ge\nu(2,l)\quad\text{for all}\quad i\ge0,\
0\le l\le s\ .\label{e36}\end{equation} This excludes possible
vertices $v=(0,j)$ or $(1,j)$ with $0<j<s$ of the trapezes, since,
otherwise, by formula (\ref{e31ms}) and relation (\ref{e37}), we
would get $\nu(v)=\nu_v\ge\min\{\mu_{0,j-1},\mu_{1,j-1}\}$,
contradicting (\ref{e36}) and the strict decrease of $\nu$ with
respect to the second variable.

Next we observe that the points $(0,0)$ and $(1,0)$ cannot serve as
vertices of the lowest trapeze. Indeed, the values of $\mu_{00}$ and
$\mu_{10}$ come from the second summand in (\ref{e31ms}), and hence
are equal to $\mu$. Thus, if the lowest trapeze is not a rectangle,
its lower side has vertex $(0,0)$. So, suppose that $(0,0)$ is a
vertex of the lower trapeze and bring this assumption to
contradiction. The upper left vertex of that lowest trapeze must be
$(2,j)$ with some $1\le j\le s$. Thus, we have
$$0\le\nu_{2j}=\nu(2,j)\le\mu\frac{s-j}{s}\ ,$$ and by linearity of
$\nu$ in each trapeze,
$$\nu(i,l)=\nu_{2j}\frac{l}{j}+\mu\frac{j-l}{j}\quad\text{for all}\quad
i\ge 2,\ 0\le l\le j\ .$$ Combining this bound with (\ref{e36}), we
derive that
$$\nu_{il}\ge\nu_{2j}\frac{l-1}{j}+\mu\frac{j-l+1}{l}>\nu_{2j}\frac{l}{j+1}+\mu\frac{j+1-l}{j+1},\quad\text{for all}\ i=0,1,\ 1\le l\le
j\ .$$ On the other hand, from (\ref{e31ms}), one obtains
$\nu(0,j+1)=\mu_{0j}=\nu_{2j}$, and hence the preceding bound
together with the convexity of $\nu$ will imply that
$$\nu_{il}>\nu(i,l)\quad\text{for all}\quad i=0,1,\ 0<l\le j\ .$$
In particular, we obtain that the subdivision of $\del_2$ contains
the triangle $\quad$ $\quad$ $\quad$ $\quad$ $\quad$
\mbox{$\del=\conv\{(0,0),(2,j),(0,j+1)\}$} (see Figure \ref{f2}(d)),
and, moreover, the limit polynomial $f_\del$ contains only three
monomials corresponding to the vertices of $\del$. It is easy to see
that then the limit curve $K_\del$ is of geometric genus
$\#(\Int(\del)\cap\Z^2)>0$, which implies that the curve $C_t$ has
handles in a neighborhood of $L'$ contradicting Lemma \ref{l1x}(vi).

\smallskip

{\it (3C)} Let us show that there is only one rectangle with
horizontal/vertical edges in the subdivision of $\del_2$. Suppose
that there are several rectangles like this (cf. Figure
\ref{f2}(e)). Let $(2,j)$, $0<j<s$, be the left lower vertex of the
upper rectangle. Arguing as in the preceding paragraph, one can
derive that
\begin{eqnarray}&\nu_{1l}>\nu(1,l)\quad\text{for all}\quad 0\le l\le
s\ ,\label{e38}\\
&\nu_{0,j+1}=\nu(0,j+1)=\nu_{2,j}=\nu(2,j)=\mu_{0,j}\ ,\nonumber\\
&\nu_{0l}\ge\nu(0,l)=\nu(2,l-1)\quad\text{for all}\quad j<l\le s\
.\nonumber\end{eqnarray} From this we conclude that the
parallelogram $\del_0=\conv\{(0,j+1),(0,s+1),(2,j),(2,s)\}$ is a
part of the $\nu$-subdivision, and that the coefficients of
$f_{\del_0}$ are nonzero only along the vertical sides of $\del_0$
and they all come from the coefficients of $t^{\nu(0,l)}$ in the
coefficients of $y^{l-1}$ in $\widetilde G_t(x,y)$, $j<l\le s$.
Particularly, $f_{\del_0}=(\sum_{i=0}^{s-j}c_iy^i)(y+x^2)$, and
hence the curve $K_{\del_0}$ splits off the conic $S_0=\{y+x^2=0\}$
and $s-j$ horizontal lines, each crossing $S_0$ transversally at two
distinct points different from the origin. By a suitable variable
change $(x,y)\mapsto (xt^{\lambda_1},yt^{\lambda_2})$ in $F_t(x,y)$
and division by the minimal power of $t$ in the obtained polynomial
(still denoted by $F_t$), we can make the corresponding function
$\nu$ vanishing along $\del_0$, which means that
$F_t(x,t)=f_{\del_0}+O(t)$. At the same time this operation takes
the conic $S$ into a family of conics $S_t=\{y+x^2+xy\cdot
O(t)=0\}$. Notice that each line of $K_{\del_0}$ crosses $S_t$ at
two points convergent to its intersection points with $S_0$. The
statement of Lemma \ref{l1x}(vi) yields that the curves $C_t$,
$t\ne0$, do not intersect $E$ in a neighborhood of $L'$. Hence, the
intersection points of the lines of $K_{\del_0}$ with $S_0$ must
smooth out in the deformation $\bigcup_\del K_\del\to K_t$, $t\ne0$.
However such a smoothing will develop at least one handle of $C_t$
in a neighborhood of $L'$, which is a contradiction to Lemma
\ref{l1x}(vi).

\smallskip

{\it (3D)} To complete the proof of the claim that the
$\nu$-subdivision of $\Delta$ is the one shown in Figure
\ref{f2}(b), it remains to exclude the subdivision depicted in
Figure \ref{f2}(f), and this can be done in the same manner as in
the preceding paragraph, where we excluded from $\nu$-subdivision
parallelograms with vertices $(0,j+1),(0,s+1),(2,j),(2,s)$.

\smallskip

{\it Step 4: Limit curves.} As we observed above, the limit curve
corresponding to the rectangle splits into the union of vertical and
horizontal lines. For the limit polynomial $f_\del$ of the trapeze
$\del=\conv\{(0,0),(0,s+1),(2,0),(2,s)\}$, the consideration of part
(3) in Step 3, notably, (\ref{e38}), gives that $f_\del(x,y)$ does
not contain monomials $xy^l$, and the coefficients of the monomials
$y^{l+1}$ and $x^2y^l$ equal the coefficient of $t^{\nu(0,l+1)}$
standing at the monomial $y^l$ in $\widetilde G_t(x,y)$ for all
$l=0,...,s-1$. Then
\begin{equation}f_\del(x,y)=c+(y+x^2)f(y)\ ,\label{e40}\end{equation}
where $c$ is defined by (\ref{e32}), and
$f(y)=y^s+b_{s-1}y^{s-1}+...+b_0$. We call this polynomial the {\it
deformation pattern} for the component $s'L'$ of $C$ associated with
the parameterized branch $V$ of $V(\overline\bp,C)$.

\begin{lemma}\label{l3}
(1) Given a curve $C=E\cup\widetilde C$ and a set $\bz'\subset
\widetilde C\cap E$ as chosen in section \ref{sec253}(3), the
coefficient $c$ is determined uniquely.

(2) For any fixed $c\ne0$, the set ${\mathcal P}(L')$ of the
polynomials given by (\ref{e40}) and defining a rational curve
consists of $s+1$ elements, and they all can be found from the
relation
\begin{equation}yf(y)+c=
\frac{c}{2}\left(\cheb_{s+1}\left(\frac{y}{(2^{s-1}c)^{1/(s+1)}}+y'\right)+1\right)
\ ,\label{e41}\end{equation} where $\cheb_{s+1}(y)=\cos((s+1)\arccos
y)$, and $y'$ is the only positive simple root of $\quad$
\mbox{$\cheb_{s+1}(y)-1$}.

(3) The germ at $f_\del\in{\mathcal P}(L')$ of the set of those
polynomials with Newton quadrangle
$\conv\{(0,0),(0,s+1),(2,0),(2,s)\}$ which define a rational curve,
is smooth of codimension $s$, and intersects transversally (at one
point $\{f_\del\}$) with the space of polynomials
$$\left\{(y+x^2)\left(y^s+\sum_{l=0}^{s-1}c_ly^l\right)\ :\
c_0,...,c_{s-1}\in\C\right\}\ .$$
\end{lemma}

{\bf Proof.} The first claim follows from the construction.

For the second claim, notice that $f_\del(-x,y)=f_\del(x,y)$ (cf.
(\ref{e40})); hence the curve $K_\del$ is the double cover of a
nonsingular curve ramified at the intersection points with the toric
divisors corresponding to the vertical sides of $\del$. It is easy
to show that the double cover is rational (cf. Lemma \ref{l1x}(vi)),
if and only if $f(y)$ has precisely $[s/2]$ double roots, and
$yf(y)+c$ has precisely $[(s+1)/2]$ double roots. Hence up to linear
change in the source and in the target, all such polynomials
$yf(y)+c$ must coincide with the \emph{Chebyshev polynomial}
$\cheb_{s+1}(y)$, which is characterized by the following
properties:

- it is real, has degree $s+1$, is real, and its leading coefficient
is $2^s$,

- for even $s$, it is odd, has $s/2$ maxima on the level $1$ and
$s/2$ minima on the level $-1$.

- for odd $s$, it is even, has $(s-1)/2$ maxima on the level $1$ and
$(s+1)/2$ minima on the level $-1$.

Thus, for $yf(y)+c$ we obtain precisely $s+1$ possibilities given by
formula (\ref{e41}). The smoothness and the dimension statement in
the third claim follow, for instance, from \cite[Theorem
6.1(iii)]{GK}. For the transversality statement, we notice that the
tangent space to $M$ consists of polynomials vanishing at each of
the $s$ nodes of $f_\del=0$, whereas the tangent space to the affine
space is $\{(y+x^2)\sum_{l=0}^{s-1}c_ly^l\}$, and its nontrivial
elements cannot vanish at $s$ points outside $y+x^2=0$ with distinct
$y$-coordinates.   \proofend

\begin{remark}\label{r6}
Similarly to Remark \ref{r5}(1), we observe that the values of $\nu$
must be integer at $(0,l)$, $0\le l\le s+1$, and hence $\mu$ must be
divisible by $s+1$.
\end{remark}

\medskip

{\bf (6) Proof of Proposition \ref{pdeg1}(2iv).} In the notation of
Proposition \ref{pdeg1}, let $\bn:C^{(i)}\to L'$ is an $s$-multiple
cover, $s\ge2$. By Proposition \ref{pnov2}(ii) and Proposition
\ref{pdeg1}(2iii), $C^{(i)}$ is rational, and the cover
$\bn:C^{(i)}\to L'$ has two ramification points, one of which $z\in
L'$ differs from the tangency point of $L'$ and $E$. In particular,
if we take a tubular neighborhood $U$ of $L'$ in $\Sigma$ and the
projections $C_t\cap U\to L'$, $t\ne0$, defined by a pencil of lines
through a generic point outside $U$, then these projections will
have ramification points converging to $z$. On the other hand, it is
easy to check that the critical points of the natural projection of
the curve $K_\del\subset\C(\del)$ given by $f_\del(x,y)=0$ with
$f_\del$ as in Lemma \ref{l3} onto the toric divisor corresponding
to the lower edge of $\del$, all lie in the big torus
$(\C^*)^2\subset\C(\del)$, and hence the critical points of the
projection $C_t\cap U\to L'$ have coordinates
$((\xi+O(t))t^\lambda,(\eta+O(t))t^{2\lambda})$ with $\lambda>0$,
$\xi,\eta\in\C$, thus, for $t\to0$ they converge to the origin,
which in our setting is the tangency point of $E$ and $L'$.
Therefore, all ramification points of the projection converge to
that tangency point, a contradiction. \proofend

\medskip

{\bf (7) Deformation patterns for nonisolated singularities, II.} In
this section we construct deformation patterns for a component
$C^{(i)}=kL(p_{kj})$ of $C$ with $k\ge1$, and we use the technique
developed in part (5) of the present section.

\begin{figure}
\setlength{\unitlength}{1cm}
\begin{picture}(15,18)(0,0)
\thinlines\put(1.5,1.5){\vector(1,0){3.5}}\put(6.5,1.5){\vector(1,0){3.5}}\put(11.5,1.5){\vector(1,0){3.5}}
\put(1.5,1.5){\vector(0,1){4}}\put(6.5,1.5){\vector(0,1){4}}\put(11.5,1.5){\vector(0,1){4}}
\put(1.5,7.5){\vector(1,0){3.5}}\put(6.5,7.5){\vector(1,0){3.5}}\put(11.5,7.5){\vector(1,0){3.5}}
\put(1.5,13.5){\vector(1,0){3.5}}\put(6.5,13.5){\vector(1,0){3.5}}\put(11.5,13.5){\vector(1,0){3.5}}
\put(1.5,7.5){\vector(0,1){4}}\put(6.5,7.5){\vector(0,1){4}}\put(11.5,7.5){\vector(0,1){4}}
\put(1.5,13.5){\vector(0,1){4.5}}\put(6.5,13.5){\vector(0,1){4.5}}\put(11.5,13.5){\vector(0,1){4.5}}
\thicklines\put(1.5,3){\line(0,1){1.5}} \put(1.5,3){\line(1,0){2.5}}
\put(4,3){\line(0,1){1.5}}\put(6.5,2){\line(0,1){2.5}}\put(7,1.5){\line(1,0){2}}\put(7,1.5){\line(-1,1){0.5}}
\put(9,1.5){\line(0,1){3}}\put(11.5,2){\line(0,1){1}}\put(11.5,2){\line(1,-1){0,5}}\put(12,1.5){\line(1,0){0.5}}
\put(12.5,1.5){\line(0,1){1.5}}\put(11.5,3){\line(1,0){1}}
\put(1.5,10){\line(0,1){0.5}}\put(1.5,10){\line(1,-2){0.5}}\put(2,9){\line(1,0){2}}
\put(4,9){\line(0,1){1.5}}\put(6.5,7.5){\line(0,1){3}}\put(6.5,7.5){\line(1,0){2.5}}
\put(9,7.5){\line(0,1){3}}\put(11.5,7.5){\line(0,1){3}}\put(11.5,7.5){\line(1,0){2.5}}
\put(14,7.5){\line(0,1){3}}
\put(1.5,16.5){\line(0,1){0.5}}\put(1.5,16.5){\line(1,-3){0.5}}
\put(2,15){\line(1,0){2}}\put(4,15){\line(0,1){2}}\put(6.5,15){\line(0,1){2}}
\put(6.5,15){\line(1,-3){0.5}}\put(7,13.5){\line(1,0){2}}\put(9,13.5){\line(0,1){3.5}}
\put(11.5,15){\line(0,1){2}}\put(11.5,15){\line(1,-3){0.5}}\put(12,13.5){\line(1,0){2}}
\put(14,13.5){\line(0,1){3.5}} \thinlines
\dashline{0.2}(4,4.5)(3.5,5)\dashline{0.2}(9,4.5)(8.5,5)
\dashline{0.2}(6.5,3)(9,3)
\dashline{0.2}(1.5,10.5)(2,11)\dashline{0.2}(6.5,10.5)(7,11)\dashline{0.2}(11.5,10.5)(12,11)
\dashline{0.2}(1.5,17)(2,17.5)\dashline{0.2}(6.5,17)(7,17.5)\dashline{0.2}(11.5,17)(12,17.5)
\dashline{0.2}(4,10.5)(3.5,11)\dashline{0.2}(9,10.5)(8.5,11)\dashline{0.2}(14,10.5)(13.5,11)
\dashline{0.2}(4,17)(3.5,17.5)\dashline{0.2}(9,17)(8.5,17.5)\dashline{0.2}(14,17)(13.5,17.5)
\dashline{0.2}(11.5,15)(12,15)\dashline{0.2}(12,15)(12.5,13.5)\dashline{0.2}(6.5,16.5)(7,15)
\dashline{0.2}(7,15)(9,15)\dashline{0.2}(11.5,16.5)(12,15)\dashline{0.2}(12,15)(14,15)
\dashline{0.2}(6.5,10)(7,9)\dashline{0.2}(7,9)(9,9)\dashline{0.2}(11.5,10)(12,9)
\dashline{0.2}(11.5,7.5)(12,9)\dashline{0.2}(12,9)(14,9)
\dottedline{0.1}(4,3)(4,1.5)
\dottedline{0.1}(1.5,9)(2,9)\dottedline{0.1}(2,9)(2,7.5)\dottedline{0.1}(4,9)(4,7.5)
\dottedline{0.1}(1.5,15)(2,15)\dottedline{0.1}(2,15)(2,13.5)
\dottedline{0.1}(4,13.5)(4,15) \put(2.8,6.5){\rm
(d)}\put(2.8,12.5){\rm (a)} \put(7.8,6.5){\rm (e)}\put(7.8,12.5){\rm
(b)}\put(12.8,6.5){\rm (f)}\put(12.8,12.5){\rm (c)}\put(2.8,6.5){\rm
(d)}\put(2.8,0.5){\rm (g)}\put(7.8,0.5){\rm (h)}\put(12.8,0.5){\rm
(i)}\put(2.6,10){$\Delta$}\put(2.6,4){$\Delta$}\put(7.6,4){$\Delta$}\put(7.6,10){$\Delta$}
\put(12.6,10){$\Delta$}\put(2.8,6.5){\rm
(d)}\put(2.8,16){$\Delta$}\put(7.8,16){$\Delta$}
\put(12.8,16){$\Delta$}\put(1.95,7.1){$1$}
\put(1.95,13.1){$1$}\put(6.95,1.1){$1$}\put(11.95,1.1){$1$}\put(12.45,1.1){$2$}\put(6.95,13.1){$1$}
\put(11.95,13.1){$1$}\put(12.45,13.1){$2$}
\put(1.1,8.9){$k$}\put(1.1,2.9){$k$}\put(6.1,2.9){$k$}\put(11.1,2.9){$k$}\put(6.2,1.9){$1$}\put(11.2,1.9){$1$}
\put(1.1,14.9){$k$}\put(6.1,14.9){$k$}\put(11.1,14.9){$k$}
\put(0.4,9.85){$2+k$}\put(5.4,9.85){$2+k$}\put(10.4,9.85){$2+k$}
\put(0.9,16.35){$2k$}\put(5.9,16.35){$2k$}\put(10.9,16.35){$2k$}
\put(3.5,7.1){$d-k$}\put(3.5,1.1){$d-k$}\put(8.5,1.1){$d-k$}\put(8.5,7.1){$d-k$}\put(13.5,7.1){$d-k$}
\put(3.5,13.1){$d-k$}\put(8.5,13.1){$d-k$}\put(13.5,13.1){$d-k$}
\end{picture}
\caption{Deformation patterns for $kL(p_{kl})$ and
$kL(z)$}\label{f3}
\end{figure}
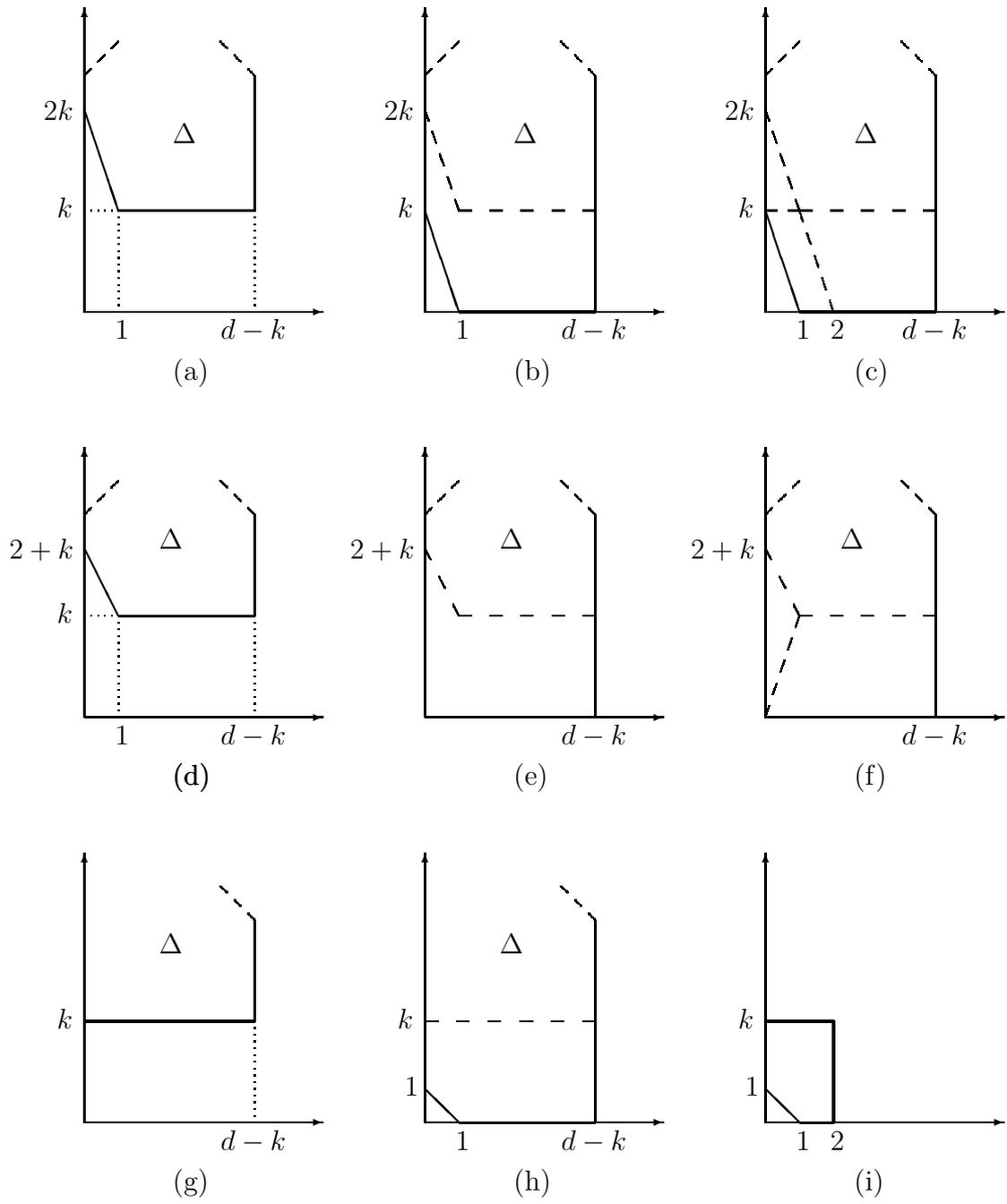

Perform the blow-down $\Pi:\Sig\to\PP^2$, contracting
$E_1,...,E_{a+1}$, and take an affine plane $\C^2\subset\PP^2$ with
the coordinates $(x,y)$ such that $\Pi(L(p_{kl}))=\{y=0\}$,
$\Pi(E_1)$ is the infinitely far point of $\Pi(L(p_{kl}))$, and the
conic $\Pi(E)$ crosses $\Pi(L(p_{kl}))$ transversally at
$\Pi(p_{kl})=(0,0)$ and at $(1,0)$, the latter point being smoothed
out in the deformation of $C$ along $V$. There exists a unique
transformation of $\C^2$
$$(x,\ y)\mapsto\left(x+\sum_{1\le i\le k}c_iy^i,\ y\right)\ ,$$ such that the
equation of $\Pi(E)\cap\C^2$ becomes
$$S(x,y):=y^k+x-x^2+\sum_{j>k}c_{0j}y^j+\sum_{j>0}(c_{1j}xy^j+c_{2j}x^2y^j)=0\
.$$ Thus, the curve $\Pi(C)\cap\C^2$ is defined by a polynomial of
the form $F_0(x,y)=S(x,y)\widetilde G_0(x,y)$ with Newton polygon
$\Delta$, whose lower part is the union of the segments
$[(0,2k),(1,k)]$ and $[(1,k),(d-k,k)]$ (shown by fat line in Figure
\ref{f3}(a)). Here the root of the truncation of $F_0$ on the
segment $[(0,2k),(1,k)]$ corresponds to the branch of $\Pi(E)$
centered at the origin, and the roots of the truncation of $F_0$ on
the segment $[(1,k),(d-k,k)]$ correspond to the intersection points
of $\Pi(E)$ with $\Pi(L(p_{kl}))$ (equal to $(1,0)$) and the
intersection points of $\Pi(L(p_{kl}))$ with the other components of
$\Pi(C)$. As in Step 1 of part (5), the curves $\Pi(C_t)\cap\C^2$
are defined by polynomials $F_t(x,y)$ converging to $F_0(x,y)$ as
$t\to0$ and given by (\ref{e31ms}) with the Newton polygon $\Delta'$
being the convex hull of the union of $\Delta$ with the points
$(0,k)$, $(1,0)]$, and $(d-k,0)$ (see Figure \ref{f3}(b)).

Then, using Lemma \ref{l1x}(i,iii,v) and proceeding as in Step 3 of
part (5), we can derive that the tropical limit of the family
$F_t(x,y)$ is as follows: \begin{itemize}\item the area
$\Delta'\setminus\Delta$ is subdivided as depicted in Figure
\ref{f3}(c),
\item the limit curve $K_{\del_1}$ with the Newton parallelogram $\del_1$
splits into two components, {\it i.e.} its defining polynomial is
$(x-b')(y^s-c'x)$ with some $b',c'\ne0$, \item the limit curve
$K_{\del_2}$ with the Newton trapeze $\del_2$ is defined by the
polynomial $\quad$ \mbox{$((x-1)y^s+c'x)\psi(x)$}, where the roots
of $\psi(x)$ correspond to the intersection points (counting
multiplicities) of $\Pi(L(p_{kl}))$ with the other components of
$\Pi(C)$.
\end{itemize} The convex piece-wise linear function
$\nu:\Delta'\to\R$ defining the above subdivision (cf. Step 2 of
part (5)) is characterized by its values
$$\nu\big|_\Delta=0,\quad\nu(0,k)=\xi>0,\quad\nu(2,0)=\eta>0\ .$$ Particularly,
this yields that $\nu(1,0)=\xi+\eta$ and $\nu(u,0)=\eta$ as $u\ge
2$. Combining this with (\ref{e31ms}) and the fact that the curves
$C_t$, $t\ne0$, cross $E$ at $p_{kl}$ and in a neighborhood of some
intersection points of $E$ with the components of $C$ different from
$L(p_{kl})$, we derive that $\xi=\eta=\mu$ with $\mu$ defined by
(\ref{e22}), and that the limit curves are as follows:
\begin{enumerate}\item[(i)] $\{y^k+x-b'=0\}$
for the triangle
$\conv\{(0,k),(0,2k),(1,k)\}$,
\item[(ii)] $\{(x-b')(y^k-c'x)=0\}$ for the parallelogram
$\del_1=\conv\{(0,k),(1,0),(2,0),(1,k)\}$,
\item[(iii)] $\{((x-1)y^k+c'x)\psi(x)=0\}$ for the trapeze $\del_2=\conv\{(2,0),(1,k),(d-k,0),(d-k,k)\}$. \end{enumerate}
The polynomial $(x-b')(y^k-c'x)$ (which uniquely determines the
polynomials in (i) and (iii)) is called the {\it deformation
pattern} for the component $kL(p_{kl})$ of $C$ associated with the
parameterized branch $V$ of $V(\overline\bp,C)$.

\begin{lemma}\label{l49}
The parameters $b'$ and $c'$ in the polynomials (i)-(iii) are
uniquely (up to a common factor) defined by the curve $C$, by the
choice of the set $\bz'$ (cf. part (5)) and by relation (\ref{e22}).
\end{lemma}

{\bf Proof.} Coming back to the coordinates $(u,v)$ introduces in
part (4), we can see that equation (\ref{e21}) and relation
(\ref{e22}) yield linear equations for $b'$ and $c'$, which can be
uniquely resolved. \proofend

\medskip

{\bf (8) Deformation patterns for nonisolated singularities, III.}
In this section we construct deformation patterns for a component
$C^{(i)}=kL(z)$ of $C$ with $k\ge1$ and $\gam^{(i)}=e_k$, and again
we apply the technique of parts (5) and (7).

Perform the blow down $\Pi:\Sig\to\PP^2$, contracting
$E_1,...,E_{a+1}$, and take an affine plane $\C^2\subset \PP^2$ with
coordinates $x,y$ such that $\Pi(L(z))\cap\C^2=\{y=0\}$, $\Pi(E_1)$
is the infinitely far point of $\Pi(L(z))$. Furthermore, we assume
that the intersection point of $\Pi(E)$ and $\Pi(L(z))$ belonging to
$\Pi(\bz')$ has coordinates $(0,0)$, the other intersection point of
$\Pi(E)$ and $\Pi(L(z))$ has coordinates $(1,0)$, and the point
$\Pi(z)\in\Pi(L_j)$ has coordinates $(x_0,0)$, $x_0\ne0,1$. At last
$\Pi(E)$ is a conic, whose equation can be made $y^2-x+x^2=0$ while
keeping the preceding data.

The curve $\Pi(C)\cap\C^2$ is then defined by a polynomial of the
form $F_0(x,y)=S(x,y)\widetilde G_0(x,y)$ with Newton polygon
$\Delta$, whose lower part is the union of the segments
$[(0,2+k),(1,k)]$ and $[(1,k),(d-k,k)]$ (shown by fat line in Figure
\ref{f3}(d)). Here the root of the truncation of $F_0$ on the
segment $[(0,2+k),(1,k)]$ corresponds to the branch of $\Pi(E)$
centered at the origin. In turn, the truncation of $F_0$ on the
segment $[(1,k),(d-k,k)]$ equals $y^kx(x-1)f(x)$, where the roots of
$f(x)$ correspond to the intersection points of $\Pi(L(z))$ with the
other components of $\Pi(C\setminus E)$. As in Step 1 of part (5),
the curves $\Pi(C_t)\cap\C^2$ are defined by polynomials $F_t(x,y)$
converging to $F_0(x,y)$ as $t\to0$ and given by (\ref{e31ms}) with
the Newton polygon $\Delta'=\conv(\Delta\cup\{(0,0),(d-k,0)\})$ (see
Figure \ref{f3}(e)).

Then, using Lemma \ref{l1}(2i,iv,v) and proceeding as in Step 3 of
part (5), we can derive that the tropical limit of the family
$F_t(x,y)$ is as follows: \begin{enumerate}\item[(i)] the area
$\Delta'\setminus\Delta$ is subdivided as depicted in Figure
\ref{f3}(f),
\item[(ii)] the limit curve $K_{\del_1}$ with the Newton trapeze
\begin{equation}\del_1=\conv\{(0,0),(d-k,0),(1,k),(d-k,k)\}\label{edel1}\end{equation} is defined by a
polynomial $F_{\del_1}(x,y)=(x-1)f(x)h^{(1)}_z(x,y)$ (with $f(x)$
introduced above), $h^{(1)}_z(x,y)$ has Newton triangle
$\conv\{(0,0),(1,0),(1,k)\}$, the curve $\{h^{(1)}_z(x,y)=0\}$
crosses $\Pi(L(z))$ at $\Pi(z)=(x_0,0)$ and crosses the line $x=1$
at one point with multiplicity $k$;
\item[(iii)] the convex piece-wise linear function
$\nu:\Delta\cup\del_1\cup\del_2\to\R$ whose graph projects to the
above subdivision (cf. Step 2 of part (5)) is characterized by its
values
$$\nu\big|_\Delta=0,\quad\nu(u,0)=\mu,\ 0\le u\le d-k\ .$$
\end{enumerate} Moreover, using Lemma \ref{ms2}, we can easily derive
that \begin{enumerate}\item[(iv)] the limit curve $K_{\del_2}$ with
the Newton triangle
\begin{equation}\del_2=\conv\{(0,0),(0,2+k),(1,k)\}\label{edel2}\end{equation} is defined by a
trinomial.\end{enumerate}

We call the polynomial $h^{(1)}_z$ in (ii) the {\it deformation
pattern} for the component $kL(z)$ of $C$ associated with the
parameterized branch $V$ of $V(\overline\bp,C)$.

\begin{lemma}\label{llj}
(1) The limit curve $K_{\del_2}$, the polynomial $f(x)$ and the
coefficients $c_{00}^{z}$, $c_{10}^{z}$, $c_{1k}^{z}$ of $1$, $x$,
and $xy^k$ in $h^{(1)}_z(x,y)$ are uniquely defined by the curve
$C$, by the choice of the set $\bz'$, and by relation (\ref{e22}).
Within these data, there are precisely $k$ polynomials
$h^{(1)}_z(x,y)$ meeting conditions (ii) above.

(2) In the space of polynomials with Newton trapeze $\del_1$ (cf.
(\ref{edel1})), the germ at $F_{\del_1}$ of the family of
polynomials, splitting into a product of binomials and a polynomial
with Newton triangle $\del_2$ (cf. (\ref{edel2})), which meets
condition (ii) above, is smooth, and it intersect with the affine
space of polynomials, having the same coefficients on the segments
$[(1,k),(d-k,k)]$ and $[(0,0),(1,k)]$ as in $F_{\del_1}$ and
vanishing at the point $(x_0,0)$, transversally at one element
(equal to $F_{\del_1}$).
\end{lemma}

{\bf Proof.} Straightforward. \proofend

We denote the set of polynomials $h^{(1)}_z$ as in Lemma \ref{llj}
by ${\mathcal P}^{(1)}(L(z))$.

\begin{remark}\label{r7}
Observe that none of the coefficients of the polynomials
$h^{(1)}_z(x,y)$ in Lemma \ref{llj} vanishes, and hence the values
$\nu(1,i)=\mu(k-i)/k$, $0\le i\le k$, must be integral, in
particular, $k$ divides $\mu$.
\end{remark}

\medskip

{\bf (9) Deformation patterns for nonisolated singularities, IV.} In
this section, using the above approach, we construct deformation
patterns for a component $C^{(i)}=kL(z)$ of $C$ with $k\ge1$ and
$\gam^{(i)}=2e_k$.

Perform the blow down $\Pi:\Sig\to\PP^2$, contracting
$E_1,...,E_{a+1}$, and take an affine plane $\C^2\subset \PP^2$ with
coordinates $x,y$ such that $\Pi(L(z))\cap\C^2=\{y=0\}$,
$\Pi(z)=(0,0)$, $\Pi(E_1)$ is the infinitely far point of
$\Pi(L(z))$, $\Pi(L')$ is the infinitely far line, and its tangency
point with $\Pi(E)$ is the infinitely far point of the axis
$\{x=0\}$. In addition to these assumptions, the equation of
$\Pi(E)$ can be brought to the form $S(x,y):=y+P(x)=0$ with a
nondegenerate quadratic polynomial $P(x)$. The curve
$\Pi(C)\cap\C^2$ is then defined by a polynomial of the form
$F_0(x,y)=S(x,y)\widetilde G_0(x,y)$ with Newton polygon $\Delta$,
whose lower part is the segment $[(0,k),(d-k,k)]$ (see Figure
\ref{f3}(g)). Here the root of the truncation of $F_0$ on the
segment $[(0,k),(d-k,k)]$ equals $y^kP(x)f(x)$, where the roots of
$f(x)$ correspond to the intersection points of $\Pi(L(z))$ with the
other components of $\Pi(C\setminus E)$. As in Step 1 of part (5),
the curves $\Pi(C_t)\cap\C^2$ are defined by polynomials $F_t(x,y)$
converging to $F_0(x,y)$ as $t\to0$ and given by (\ref{e31ms}) with
the Newton polygon
$\Delta'=\conv(\Delta\cup\{(1,0),(0,1),(d-k,0)\})$ (see Figure
\ref{f3}(h)). Using Lemma \ref{l1x}(i,v) and proceeding as in Step 3
of part (5), we can derive that the tropical limit of the family
$F_t(x,y)$ is as follows:
\begin{enumerate}\item[(i)] $\Delta'$ is
subdivided as depicted in Figure \ref{f3}(h);
\item[(ii)] the limit curve $K_{\del'}$ with the Newton polygon
$\del'=\Delta'\setminus\Delta$ is defined by a polynomial
$F_{\del'}(x,y)=f(x)h_z^{(2)}(x,y)$ ($f(x)$ introduced above),
$h^{(2)}_z(x,y)$ has Newton polygon
$\conv\{(0,1),(1,0),(2,0),(2,k),(0,k)\}$;
\item[(iii)] the convex piece-wise linear function
$\nu:\Delta'\to\R$, whose graph projects to the subdivision in
Figure \ref{f3}(h), is characterized by its values
$$\nu\big|_\Delta=0,\quad\nu(u,0)=\mu,\ 1\le u\le d-k\ .$$
\end{enumerate} Using Lemma \ref{ms2}, we can easily derive that
\begin{equation}h^{(2)}_z(x,y)=y\varphi(y)P(x)+Q(x),\quad
\varphi(y)=y^{k-1}+\sum_{j=0}^{k-2}c_jy^j\ .\label{eG}\end{equation}
Finally, we observe that by Lemma \ref{l1}(2v),
$\{h^{(2)}_z(x,y)=0\}$ is a rational curve. We call the polynomial
$h^{(2)}_z$ the {\it deformation pattern} for the component $kL(z)$
of $C$ associated with the parameterized branch $V$ of
$V(\overline\bp,C)$.

\begin{lemma}\label{llz}
(1) The coefficients of the polynomial $Q(x)=b_1x+b_2x^2$ in
(\ref{eG}) are uniquely defined by the curve $C$, by the choice of
the set $\bz'$, and by relation (\ref{e22}). Within these data,
there are precisely $k^2$ polynomials $h^{(2)}_z(x,y)$ satisfying
(\ref{eG}) and defining a rational curve.

(2) In the space of polynomials with Newton polygon $\del'$ (see
Figure \ref{f3}(h)), the germ at $F_{\del'}$ of the family of
polynomials, splitting into a product of binomials and a polynomial
with Newton polygon $\del$, which defines a rational curve, is
smooth, and it intersects transversally at one element (equal to
$F_{\del'}$) with the affine space of polynomials $g(x,y)f(x)$,
where $g(x,y)$ is given by formula (\ref{eG}) with free parameters
$c_0,...,c_{k-2}$.
\end{lemma}

{\bf Proof.} We explain only the number $k^2$ of polynomials
$h^{(2)}_z(x,y)$ in count. Computing singularities of
$\{h^{(2)}_z(x,y)=0\}$ via relation (\ref{eG}) and equations
$$h^{(2)}_z(x,y)=\frac{\partial}{\partial x}h^{(2)}_z(x,y)=\frac{\partial}{\partial y}h^{(2)}_z(x,y)=0\ ,$$
we, first, obtain an equation $QP'-Q'P=0$, which has two distinct
roots $x_1,x_2$ (it follows from the fact that, by construction, the
root of $Q(x)$ can be made generic with respect to the roots of
$P(x)$), and, second, we derive that the polynomial $y\varphi(y)$ in
(\ref{eG}) has precisely two critical levels (equal to
$-Q(x_1)/P(x_1)$ and $-Q(x_2)/P(x_2)$), and hence (cf. the proof of
Lemma \ref{l3}), up to a uniquely defined affine modification in the
source and in the target, and up to a shift, the polynomial
$y\varphi(y)$ coincides with one of the $k$ Chebyshev polynomials of
degree $k$. At last, the vanishing at $y=0$ determines precisely $k$
possible shifts for each of the above modified Chebyshev
polynomials. \proofend

We denote the set of polynomials $h^{(2)}_z$ as in Lemma \ref{llj}
by ${\mathcal P}^{(2)}(L(z))$.

\begin{remark}\label{r8}
Observe that none of the coefficients of the polynomials
$h^{(2)}_z(x,y)$ in Lemma \ref{llj} vanishes, and hence the values
$\nu(1,i)=\mu(k-i)/k$, $0\le i\le k$, must be integral, in
particular, $k$ divides $\mu$.
\end{remark}

\medskip

{\bf (10) From deformation patterns to deformation.} Observe that
formula (\ref{e18}) counts the number of possible collections of
deformation patterns as described in Lemmas \ref{l2}, \ref{l3},
\ref{l49}, \ref{llj}, and \ref{llz} for all choices made in part (3)
of the present section. Thus we complete the proof of Proposition
\ref{p5}(2) with the following

\begin{lemma}\label{l4}
In the notations and hypotheses of Proposition \ref{p5}(2), let us
be given the set $\bz'$, chosen as in part (3). Then

(1) for any collection of polynomials
\begin{equation}\begin{cases}&\Psi_{ij}(x,y)\in{\mathcal
P}(q'_{ij})\quad \text{for all}\quad q'_{ij}\in\bz',\ \text{as in
Lemma \ref{l2}},\\ &\Psi'(x,y)\in{\mathcal P}(L')\ \text{and}\
\Psi''(x,y)\in{\mathcal P}(L''),\
\text{as in Lemma \ref{l3}},\\
&h_z^{(1)}\in{\mathcal P}^{(1)}(L(z))\quad\text{for all}\
C^{(i)}=kL(z)\ \text{with}\ \gam^{(i)}=e_k,\ \text{as in Lemma
\ref{llj}},\\ &h_z^{(2)}\in{\mathcal P}^{(2)}(L(z))\quad\text{for
all}\ C^{(i)}=kL(z)\ \text{with}\ \gam^{(i)}=2e_k,\ \text{as in
Lemma \ref{llz}},\end{cases}\label{e42}\end{equation} there exists a
unique normally parameterized branch $V$ of $V(\overline\bp,C)$
realizing the given polynomials as deformation patterns;

(2) there is one-to-one correspondence between the curves $C'\in
V(\overline\bp,C)$ passing through a given point $p'\in L_p$ and the
collections of polynomials (\ref{e42}).
\end{lemma}

{\bf Proof.} By technical reason, we shall work with
$H^0(\Sig,{\mathcal O}_\Sig(D))$ rather than with the linear system
$|D|$, and the curves $C'\in|D|$ will be lifted to sections in
$H^0(\Sig,{\mathcal O}_\Sig(D))$ denoted by $F_{C'}$ and specified
below.

Using the data of Lemma \ref{l4} and given deformation patterns
(\ref{e42}), we will derive a provisional formula for the lift of
$V$ to $H^0(\Sig,{\mathcal O}_\Sig(D))$. Then, using certain
transversality conditions, we will show that the formula does define
a unique parameterized branch $V$ of $V(\overline\bp,C)$.

\smallskip

{\it Step 1.} In the linear system $|D|$, consider the germ at $C$
of the family of curves, which \begin{itemize}\item pass through
$\overline\bp\setminus\{z\ :\ L(z)\subset C\}$,
\item in neighborhood of the points of
$$\bp\cup(\Sing(C_{\redu})\setminus
E)\cup(\bz\setminus\{q_{kl}\in\bz\ :\ q_{kl}\in L(z),\ L(z)\subset
C\})\ ,$$ realize local deformations as described in Lemma
\ref{l1}(2i,ii,iii,iv),
\item in neighborhood of the points $\{q_{kl}\in\bz\ :\ q_{kl}\in L(z),\ L(z)\subset C\}$, realize
deformations in which the local branches of $E$ do not glue up with
local branches of the lines $L(z)$ (weaker that that in Lemma
\ref{l1x}(iv), where we additionally require a special tangency
condition).
\end{itemize} We will show that this family lifts to a smooth
variety germ $M\subset H^0(\Sig,{\mathcal O}_\Sig(D))$. Namely, we
will describe $M$ as an intersection of certain smooth germs and
then verify the transversality of the intersection. In the following
computations we use the model of part (4) with affine coordinates
$u,v$.

In a neighborhood of a point
$p_{kl}\in\bp\setminus\bigcup_iC^{(i)}$, in the coordinates
$x=u-u(p_{kl})$, $y=v$, we have
$$F_{C}(x,y)=y(\widetilde c+\psi(x,y)),\quad \widetilde c\ne0,\ \psi(0,0)=0\ .$$ By
Lemma \ref{l1}(ii), the local deformation of $C$ along
$V(\overline\bp,C)$ can be described as
$$y(\widetilde c+\psi(x,y))+\sum_{kk'+l'\ge
k}c_{k'l'}x^{k'}y^{l'},\quad c_{k'l'}\in(\C,0)\ ,$$ which defines a
germ of a linear subvariety $M_{p_{kl}}\subset H^0(\Sig,{\mathcal
O}_\Sig(D))$ which embeds into ${\mathcal O}_{\Sig,p_{kl}}$ as the
intersection of the ideal $I_{p_{kl}}=\langle y,x^k\rangle$ with the
image of $H^0(\Sig,{\mathcal O}_\Sig(D))$ \footnote{There is no
canonical map $H^0(\Sig,{\mathcal O}_\Sig(D))\to{\mathcal
O}_{\Sig,p_{kl}}$, but the image of $H^0(\Sig,{\mathcal O}_\Sig(D))$
in ${\mathcal O}_{\Sig,p_{kl}}$ is defined correctly.}.

In a neighborhood of a point $p_{kl}\in\bp\cap C^{(i)}$, where
$C^{(i)}$ is reduced, in the coordinates $x=u-u(p_{kl})$, $y=v$, we
have
\begin{equation}F_{C}(x,y)=y(\widetilde cy+\widetilde c'x^k+\psi(x,y)),\quad
\widetilde c\widetilde c'\ne0,\
\psi(x,y)=\sum_{k'+kl'>2k}b_{k'l'}x^{k'}y^{l'}\
.\label{enn1}\end{equation} By Lemma \ref{l1}(ii), the local
deformation of $C$ along $V(\overline\bp,C)$ can be described as
$$\left(y+\sum_{k',l'\ge
0}c_{k'l'}x^{k'}y^{l'}\right)\left(\widetilde cy+\widetilde
c'x^k+\psi(x,y)+\sum_{k'+kl'\ge 2k}c'_{k'l'}x^{k'}y^{l'}\right)\ ,$$
$c_{k'l'},c'_{k'l'}\in(\C,0)$. It is immediate that this formula
defines the germ of a smooth subvariety $M_{p_{kl}}\subset
H^0(\Sig,{\mathcal O}_\Sig(D))$, whose tangent space embeds into
${\mathcal O}_{\Sig,p_{kl}}$ as the intersection of the ideal
\mbox{$I_{p_{kl}}=\langle \widetilde cy+\widetilde
c'x^k+\psi(x,y),y^2,x^ky\rangle$} with the image of
$H^0(\Sig,{\mathcal O}_\Sig(D))$. In the same way, for a point
$p_{kl}\in\bp\cap C^{(i)}$, where $C^{(i)}=kL(p_{kl})$, the
considered deformation as in Lemma \ref{l1x}(iii) is realized in a
smooth variety germ $M_{p_{kl}}\subset H^0(\Sig,{\mathcal
O}_\Sig(D))$, whose tangent space embeds into ${\mathcal
O}_{\Sig,p_{kl}}$ as the intersection of the ideal
\mbox{$I_{p_{kl}}=\langle x^k+\psi(x,y),y^2,x^ky\rangle$} with the
image of $H^0(\Sig,{\mathcal O}_\Sig(D))$.

In a neighborhood of a point $q_{kl}\in\bz\cap C^{(i)}$, where
$C^{(i)}$ is reduced, in the coordinates $x=u-u(q_{kl})$, $y=v$, we
again have (\ref{enn1}), and by Lemma \ref{l1}(iv), the local
deformation of $C$ along $V(\overline\bp,C)$ can be described as
$$\left(y+\sum_{k',l'\ge
0}c_{k'l'}x^{k'}y^{l'}\right)\left(\widetilde cy+d(x+\widetilde
c')^k+\psi(x,y)+\sum_{k'+kl'\ge 2k}c'_{k'l'}x^{k'}y^{l'}\right)\ ,$$
$c_{k'l'},c',c'_{k'l'}\in(\C,0)$. This formula also defines the germ
of a smooth subvariety $\quad\quad$ \mbox{$M_{q_{kl}}\subset
H^0(\Sig,{\mathcal O}_\Sig(D))$}, whose tangent space embeds into
${\mathcal O}_{\Sig,q_{kl}}$ as the intersection of the ideal
\mbox{$I_{q_{kl}}=\langle cy+dx^k+\psi(x,y),y^2,x^{k-1}y\rangle$}
with the image of $H^0(\Sig,{\mathcal O}_\Sig(D))$.

For a point $q_{kl}\in\bz\cap C^{(i)}$, where $C^{(i)}=kL(z)$, in
the local coordinates $x=u-u(q_{kl})$, $y=v$, where, in addition,
$L(z)=\{x=0\}$, the considered deformation can be described as
$$\left(y+\sum_{k',l'\ge
0}c_{k'l'}x^{k'}y^{l'}\right)\left(x^k+\sum_{k',l'\ge
0}c'_{k'l'}x^{k'}y^{l'}\right)\ ,$$ $c_{k'l'},c'_{k'l'}\in(\C,0)$.
This formula again defines the germ of a smooth subvariety
$\qquad\qquad$ \mbox{$M_{q_{kl}}\subset H^0(\Sig,{\mathcal
O}_\Sig(D))$}, whose tangent space embeds into ${\mathcal
O}_{\Sig,q_{kl}}$ as the intersection of the ideal
\mbox{$I_{q_{kl}}=\langle y,x^k\rangle$} with the image of
$H^0(\Sig,{\mathcal O}_\Sig(D))$.

At last, we notice that, if $z\in\Sig\setminus E$ is a center of
$l\ge2$ distinct smooth local branches $P_1,...,P_l$ of $C$ of
multiplicities $r_1,...,r_l$, respectively, and, in some local
coordinates $x,y$, we have $z=(0,0)$, $P_i=\{f_i(x,y)=0\}$,
$i=1,...,l$, then (cf. \cite[Section 5.2]{Sh0}) the closure
$M_z\subset H^0(\Sig,{\mathcal O}_\Sig(D))$ of the set of sections
defining in a neighborhood of $z$ curves with $\sum_{1\le i<j\le
l}r_ir_j(P_i\cdot P_j)$ nodes, is smooth at $F_{C}$ and its tangent
space embeds into ${\mathcal O}_{\Sig,z}$ as the intersection of the
ideal $I_z=\langle \{\prod_{i\ne j}f_i^{r_i}\ :\ j=1,...,l\}\rangle$
with the image of $H^0(\Sig,{\mathcal O}_\Sig(D))$.

The branch $V$ we are looking for must lie inside the intersection
of the above germs $M_{p_{kl}},M_{q_{kl}},M_z$ with the space of
sections vanishing at $\overline\bp\setminus\{z\ :\ L(z)\subset
C\}$. We claim that this intersection is transverse, or,
equivalently, that $H^1(\Sig,{\mathcal J}_{Z_1/\Sig}(D))=0$, where
${\mathcal J}_{Z_1/\Sig}\subset{\mathcal O}_\Sig$ is the ideal sheaf
of the scheme $Z_1\subset\Sig$ supported at
\mbox{$\bp\cup\bz\cup(\Sing(C_{\redu})\setminus E)$}, where it is
defined by the above local ideals $I_{p_{kl}},I_{q_{kl}},I_z$,
respectively, and supported at $\overline\bp\setminus\{z\ :\
L(z)\subset C\}$, where it is defined by the maximal ideals. From
the exact sequence of sheaves (cf. \cite{Hir})
$$0\to{\mathcal J}_{Z_1:E/\Sig}(D-E)\to{\mathcal J}_{Z_1/\Sig}(D)\to{\mathcal J}_{Z_1\cap E/E}(DE)\to0$$
we obtain the cohomology exact sequence
\begin{equation}H^1(\Sig,{\mathcal J}_{Z_1:E/\Sig}(D-E))\to H^1(\Sig,{\mathcal J}_{Z_1/\Sig}(D))\to H^1(E,{\mathcal J}_{Z_1\cap
E/E}(DE))\ .\label{enn3}\end{equation} Here $H^1(E,{\mathcal
J}_{Z_1\cap E/E}(D))=0$ in view of Riemann-Roch, since $\deg(Z_1\cap
E)=DE$. Thus, it remains to prove
\begin{equation}H^1(\Sig,{\mathcal J}_{Z_1:E/\Sig}(D-E))=0\
,\label{enn4}\end{equation} Notice that the scheme $Z_1:E$ coincides
with $Z_1$ in $\Sig\setminus E$, and it is defined at the points
$p_{kl}\in\bp\cap C^{(i)}$ by the ideals $J_{p_{kl}}=\langle
y,x^k\rangle\subset{\mathcal O}_{\Sig,p_{kl}}$, at the points
$q_{kl}\in\bz\setminus\{q_{kl}\ :\ q_{kl}\in L(z),\ L(z)\subset C\}$
by the ideals $J_{q_{kl}}=\langle y,x^{k-1}\rangle\subset{\mathcal
O}_{\Sig,q_{kl}}$. From the definition of the ideals $I_z$,
$z\in\Sing(\widetilde C_{\redu})\setminus E$, and the Noether
fundamental theorem, we have a canonical decomposition
$$H^0(\Sig,{\mathcal J}_{Z_1:E/\Sig}(D-E))\simeq
\bigoplus_{C^{(i)}\ne s'L',s''L''}\prod_{j\ne i}F_{C^{(j)}}\cdot
(F_{L'})^{s'}(F_{L''})^{s''}\cdot H^0(\Sig,{\mathcal
J}_{Z_1(C^{(i)})/\Sig}(C^{(i)}))$$
$$\qquad\qquad\qquad\oplus\prod_{i=1}^mF_{C^{(i)}}\cdot (F_{L''})^{s''}\cdot
H^0(\Sig,{\mathcal O}_\Sig(-s'(K_\Sig+E)))$$
\begin{equation}\qquad\qquad\qquad\;\oplus\prod_{i=1}^mF_{C^{(i)}}\cdot
(F_{L'})^{s'}\cdot H^0(\Sig,{\mathcal O}_\Sig(-s''(K_\Sig+E)))\
,\label{enn13}\end{equation} where $Z_1(C^{(i)})\subset C^{(i)}$ is
the part of $Z_1:E$ supported along $C^{(i)}$ and disjoint from the
other components of $C$ different from $E$ (in particular,
$Z_1(kL(z))=\emptyset$). It follows that (\ref{enn4}) turns into a
sequence of equalities
\begin{equation}
H^1(\Sig,{\mathcal J}_{Z_1(C^{(i)})/\Sig}(C^{(i)}))=0,\quad
i=1,...,m\ . \label{enn12}\end{equation} We intend to prove even a
stronger result: \begin{itemize} \item for a reduced
$C^{(i)}\not\in|-(K_\Sig+E)|$, we replace (\ref{enn12}) with
\begin{equation}H^1(\Sig,{\mathcal
J}_{Z_2(C^{(i)})/\Sig}(C^{(i)}))=0\ ,\label{enn8}\end{equation}
where the zero-dimensional scheme $Z_2(C^{(i)})$ is obtained from
$Z_1(C^{(i)})$ by adding the (simple) points $\overline\bp\cap
C^{(i)}$ and the scheme supported at $\bz'\cap C^{(i)}$ and defined
in $q'_{kl}\in\bz'\cap C^{(i)}$ by the ideal $J_{q'_{kl}}=\langle
y,x^{k-1}\rangle\subset{\mathcal O}_{\Sig,q'_{kl}}$; \item for
$C^{(i)}=kL(z)$ such that $\gam^{(i)}=e_k$, we replace
(\ref{enn12}), which in view of $Z_1(kL(z))=\emptyset$ reads as
$H^1(\Sig,{\mathcal O}_\Sig(kL(z)))=0$, with
\begin{equation}H^1(\Sig,{\mathcal J}_{Z^1_2(kL(z))/\Sig}(kL(z)))=0\ ,\label{enn8a}\end{equation} where the zero-dimensional
scheme $Z^1_2(kL(z))$ is supported at the point $q_{kl}\in
L(z)\cap\bz$ and is defined there by the ideal $\langle
y,x^k\rangle$; \item for $C^{(i)}=kL(z)$ such that
$\gam^{(i)}=2e_k$, we replace (\ref{enn12}) with
\begin{equation}H^1(\Sig,{\mathcal J}_{Z^2_2(kL(z))/\Sig}(kL(z)))=0\
,\label{enn8b}\end{equation} where the zero-dimensional scheme
$Z^2_2(kL(z))$ is supported at the point $z\in L(z)\cap\overline\bp$
and is defined there by the ideal $\langle x,y^k\rangle$ in the
coordinates $x,y$ introduced in part (9).
\end{itemize}

To derive (\ref{enn8}), we notice that $H^1(\Sig,{\mathcal
J}_{Z_1(C^{(i)})\setminus\overline\bp/\Sig}(C^{(i)}))=0$ is
equivalent to \begin{equation}H^1(\hat C^{(i)},{\mathcal N}_{\hat
C^{(i)}}(-\bd^{(i)}))=0\ ,\label{enn6}\end{equation} where $\hat
C^{(i)}$ is the normalization of $C^{(i)}$, and
$$\bd^{(i)}=\sum_{p_{kl}\in\bp\cap C^{(i)}}k\cdot
\hat p_{kl}+\sum_{q_{kl}\in\bz\cap C^{(i)}}(k-1)\cdot \hat
q_{kl}+\sum_{q'_{kl}\in\bz'\cap C^{(i)}}(k-1)\cdot \hat q'_{kl}$$
($\hat q'_{kl}$ being the lift of $q'_{kl}$ to $\hat C^{(i)}$). In
its turn, (\ref{enn6}) comes from (\ref{enn5}). This yields, in
particular, that (details of computation left to the reader)
$$h^0(\Sig,{\mathcal
J}_{Z_1(C^{(i)}\setminus\overline\bp/\Sig}(C^{(i)}))=h^0(\Sig,{\mathcal
O}_\Sig(C^{(i)}))-\deg(Z_1(C^{(i)})\setminus\overline\bp)=\#(\overline\bp\cap
C^{(i)})+1\ .$$ Due to the generic choice of the configuration
$\overline\bp$, we conclude that $$h^0(\Sig,{\mathcal
J}_{Z_2(C^{(i)})/\Sig}(C^{(i)}))=h^0(\Sig,{\mathcal
J}_{Z_1(C^{(i)})\setminus\overline\bp/\Sig}(C^{(i)}))-\#(\overline\bp\cap
C^{(i)})=1\ ,$$ which immediately implies (\ref{enn8}). Similarly,
relations (\ref{enn12}) for $C^{(i)}=kL(p_{kl})$ and relations
(\ref{enn8a}), (\ref{enn8b}) for $C^{(i)}=kL(z)$ lift to the
$k$-sheeted coverings $\bn:\hat L(p_{kl})\to L(p_{kl})$, $\bn:\hat
L(z)\to L(z)$ (cf. Proposition \ref{pnov2}(iii,iv)) in the form
$$H^1(\hat L(p_{kl}),{\mathcal N}_{\hat L(p_{kl})}(-\bn^*(p_{kl}))=H^1(\hat
L(z),{\mathcal N}_{\hat L(z)}(-\bn^*(q_{kl})))=H^1(\hat
L(z),{\mathcal N}_{\hat L(z)}(-\bn^*(z)))=0\ ,$$ which all come from
Riemann-Roch due to
$$\deg\bn^*(p_{kl})=\deg\bn^*(q_{kl})=\deg\bn^*(z)=k\ .$$

Thus, we conclude that the considered family $M$ is smooth.

\smallskip

{\it Step 2.} Now we construct a certain parametrization of $M$.
For, choose a special basis for its tangent space
$H^0(\Sig,{\mathcal J}_{Z_1/\Sig}(D))$.

\smallskip

{\bf (2a)} In view of (\ref{enn4}), we have
\begin{equation}H^0(\Sig,{\mathcal J}_{Z_1/\Sig}(D))\simeq F_E\cdot
H^0(\Sig,{\mathcal J}_{Z_1:E/\Sig}(D-E))\oplus H^0(E,{\mathcal
J}_{Z_1\cap E/E}(DE))\ ,\label{enn2606}\end{equation} and hence
there is a section $F\in H^0(\Sig,{\mathcal
J}_{Z_1\cup\overline\bp/\Sig}(D))$, which in the coordinates $u,v$
of part (4) can be written in the form
$$f_2(u)+v\widetilde f_2(u,v)\ ,$$ where $f_2$ is given by
(\ref{enn10}).

\smallskip

{\bf (2b)} Using (\ref{enn8}) for a reduced
$C^{(i)}\not\in|-(K_\Sig+E)|$, we get $H^0(\Sig,{\mathcal
J}_{Z_2(C^{(i)})/\Sig}(C^{(i)}))=F_{C^{(i)}}\cdot\C$, where
$F_{C^{(i)}}$ is specified so that, in the coordinates $x,y$ in a
neighborhood of each point $q'_{kl}\in\bz'\cap C^{(i)}$ (see part
(4)), the coefficient of $y^2$ in $F_{C^{(i)}}(x,y)$ equals
$c_{02}^{kl}$ as defined in part (4). Next, from the definition of
the scheme $Z_2(C^{(i)})$, we get that $H^0(\Sig,{\mathcal
J}_{Z_1(C^{(i)})/\Sig}(C^{(i)}))$, $1\le l\le m$, admits a basis
consisting of $F_{C^{(i)}}$ and the sections $F_{i,j,kl}$ labeled by
the points $q'_{kl}\in\bz'\cap C^{(i)}$ and the numbers $0\le
j<k-1$, such that \begin{itemize}\item
$\jet_{k-2,q'_{kl}}(F_{i,j,kl}(x,0))=x^j$, $0\le j<k-1$, in the
coordinates $x=u-u(q'_{kl})$, $y=v$:
\item
$\jet_{k'-2,q'_{k'l'}}(F_{i,j,kl}(x,0))=0$ in the coordinates
$x=u-u(q'_{k'l'})$, $y=v$, for all $i,j,k,l$ and
$q'_{k'l'}\in\bz'\cap C^{(i)}$, $(k',l')\ne(k,l)$.
\end{itemize}

\smallskip

{\bf (2c)} The scheme $Z_1(kL(p_{kl}))$ is determined by the
condition to cross $E$ at the point $p_{kl}$ with multiplicity $k$.
Hence $H^0(\Sig,{\mathcal
J}_{Z_1(kL(p_{kl}))/\Sig}(kL(p_{kl})))=(F_{L(p_{kl})})^k\cdot\C$.

\smallskip

{\bf (2d)} For the $k$-dimensional space $H^0(\Sig,{\mathcal
O}_\Sig(C^{(i)}))$, where $C^{(i)}=kL(z)$, we choose the basis as
follows: \begin{itemize}\item if $\gam^{(i)}=e_k$, we take
$(F_{L(z)})^j(F_{\widetilde L(z)})^{k-j}$, $0\le j\le k$, where, in
the coordinates $x,y$ of part (8), $\widetilde L(z)\in|-(K_\Sig+E)|$
is the infinitely far line (respectively, $F_{L(z)}(x,y)=y$,
$F_{\widetilde L(z)}(x,y)=1$);
\item if $\gam^{(i)}=2e_k$, we take $(F_{L(z)})^j(F_{\widetilde L(z)})^{k-j}$, $0\le j\le k$, where, in the
coordinates $x,y$ of part (9), $\widetilde L(z)\in|-(K_\Sig+E)|$ is
the infinitely far line (respectively, $F_{L(z)}(x,y)=y$,
$F_{\widetilde L(z)}(x,y)=1$).\end{itemize}

\smallskip

{\bf (2e)} Finally, for the summands $H^0(\Sig,{\mathcal
O}_\Sig(-s'(K_\Sig+E)))$ and $H^0(\Sig,{\mathcal
O}_\Sig(-s''(K_\Sig+E)))$ in (\ref{enn13}), choose the following
bases
$$F_{L'}^kF_{\widetilde L'}^{s'-j},\ 0\le j\le s',\quad\text{and}\quad
F_{L''}^kF_{\widetilde L''}^{s''-j},\ 0\le j\le s''\ ,$$ where
$F_{L'}=y$, $F_{\widetilde L'}(x,y)=1$ in the coordinates $x,y$
introduced in Step 1 of part (5), and $F_{L''}$, $F_{\widetilde
L''}$ are defined similarly.

\smallskip

Thus, our basis for $H^0(\Sig,{\mathcal J}_{z_1/\Sig}(D))$ contains
the section $F$ (see (2a)), the sections built in (2b)-(2e)
(multiplied with extra factors along (\ref{enn13}),
(\ref{enn2606})), and certain unspecified completion.

We can write then a parametrization of $M$ in the form
$$F_E\cdot\prod_{C^{(i)}\ne
kL(p_{kl}),kL(z)}\left(F_{C^{(i)}}+\sum_{ij,k}\xi_{j,k,l}F_{i,j,kl}\right)\cdot\prod_{C^{(i)}=kL(p_{kl})}(F_{L(p_{kl})})^k$$
$$\times
\prod_{C^{(i)}=kL(z)}\left((F_{L(z)})^k+\sum_{j=1}^{k}\xi_{j,z}(F_{L(z)})^{k-j}(F_{\widetilde
L(z)})^j\right)$$
$$\quad\;\times\left(F_{L'}^{s'}+
\sum_{j=0}^{s'-1}\xi'_jF_{L'}^jF_{\widetilde
L'}^{s'-j}\right)\left(F_{L''}^{s''}+
\sum_{j=0}^{s''-1}\xi''_jF_{L''}^jF_{\widetilde
L''}^{s''-j}\right)$$
\begin{equation}+\xi_0F+O(\overline\xi^2)\
,\label{enn14}\end{equation} where $\overline\xi$ is the sequence of
\begin{equation}N:=\sum_{C^{(i)}\ne kL(z)}(I\gam^{(i)}-|\gam^{(i)}|)+\sum_{C^{(i)}=kL(z)}k+s'+s''+1\label{eparam}\end{equation} free parameters
$$\xi_0,\xi_{j,s},\xi_{ijkl},\xi'_s,\xi''_s\in(\C,0)\ .$$

\smallskip

{\it Step 3.} We are ready to construct the desired parameterized
branch $V$ of $V(\overline\bp,C)$ realizing deformation patterns
(\ref{e42}). Namely, imposing additional conditions to the
parameters $\overline\xi$ in (\ref{enn14}), we will expose a formula
for a lift of the required branch $V$ to $H^0(\Sig,{\mathcal
O}(D))$.

Put \begin{equation}\mu=\lcm\left(\{k\ :\
q'_{kl}\in\bz'\}\cup\{s'+1,s''+1\}\cup\{k\ :\
kL(z)\in\{C^{(i)}\}_{i=1,...,m}\}\right)\label{efornw}\end{equation}
(cf. Remarks \ref{r5}, \ref{r6}, \ref{r7}, and \ref{r8}) and let
$\xi_0=t^\mu c(t)$ with $t\in(\C,0)$ and $c(0)=c$ given by formula
(\ref{e23}).

\smallskip

{\bf (3a)} For a point $q'_{kl}\in\bz\cap C^{(i)}$ such that
$C^{(i)}\not\in|-(K_\Sig+E)|$ is reduced, choose a deformation
pattern
$$\Psi_{kl}(x,y)=c^{kl}_{02}y^2+c^{kl}_{j1}x^jy+\sum_{j<k,\ j\equiv
k\mod2}b^{kl}_{j1}x^jy+b^{kl}_{00}$$ (cf. (\ref{e39}), (\ref{e42}),
and Lemma \ref{l2}(2)). Let $F_{kl}(x,y)$ be the section
(\ref{enn14}) written in the coordinates $x=u-u(q'_{kl})$, $y=v$.
The coefficient $c_{k-1,1}^{kl}(\overline\xi)$ of $x^{k-1}y$ in
$F_{kl}$ vanishes at $\overline\xi=0$. Since the coefficient of
$x^ky$ in $F_{kl}$ turns into $c_{k1}^{kl}\ne0$ as $\overline\xi=0$,
there is a function $\tau_{kl}(\overline\xi)$ vanishing at the
origin such that the coefficient of $x^{k-1}y$ in
$F_{kl}(x-\tau(\overline\xi),y)$ equals $0$, and hence the
coefficient of $x^jy$, $j<k-1$, in $F_{kl}(x-\tau(\overline\xi),y)$
equals $\xi_{i,j,kl}+\lambda_{i,j,kl}(\overline\xi)$ with a function
$\lambda_{i,j,kl}\in O(\overline\xi^2)$. Thus, due to the smoothness
and transversality statement in Lemma \ref{l2}(3), the conditions to
fit the given deformation pattern $\Psi_{ij}$ and to realize a local
deformation as described in Lemma \ref{l1}(v) amounts to a system of
relations
\begin{equation}\xi_{i,j,kl}+\lambda_{i,j,kl}(\overline\xi)=\begin{cases}t^{\mu(k-j)/(2k)}(b_{j1}^{kl}+O(t)),\quad
& j\equiv k\mod2,\\ O(t^{[\mu(k-j)/(2k)]+1}),\quad & j\not\equiv
k\mod2,\end{cases},\quad 0\le j<k-1\ .\label{enn15}\end{equation}

\smallskip

{\bf (3b)} Let $C^{(i)}=kL(z)$ with $\gam^{(i)}=e_k$, and let, in
the coordinates $x,y$ of part (8),
$h^{(1)}_z(x,y)=x(y^k+\sum_{j=0}^{k-1}b_{z,j}y^j)-x_0b_{z,0}\in{\mathcal
P}^{(1)}(L(z))$. In these coordinates, expression (\ref{enn14})
reads as
$$(y^2-x+x^2)\left(y^k+\sum_{j=0}^{k-1}\xi_{z,j}y^j\right)(1+O(x,y,\overline\xi))+t^\mu(x_0b_{z,0}+O(t,\overline\xi))(x-1+O(x^2,y,\overline\xi))\ .$$
In view of the smoothness and transversality statement in Lemma
\ref{llj}(2), the conditions to fit the deformation pattern
$h^{(1)}_z(x,y)$, to realize a local deformation as described in
Lemma \ref{l1}(2iv), and to hit the point $z\in
L(z)\cap\overline\bp$ altogether amount to a system of equations
\begin{equation}\xi_{z,j}+O(\overline\xi^2)=t^{\mu(k-j)/k}(b_{z,j}+O(t,\overline\xi)),\quad
0\le j<k\ .\label{elj}\end{equation} Similarly, let $C^{(i)}=kL(z)$
with $\gam^{(i)}=2e_k$, and let, in the coordinates $x,y$ from part
(9),
$h^{(2)}_z(x,y)=P(x)(y^k+\sum_{j=1}^{k-1}c_{z,j}y^j)+Q(x)\in{\mathcal
P}^{(2)}(L(z))$ be a deformation pattern (cf. (\ref{eG})). In these
coordinates, expression (\ref{enn14}) reads as
$$(y^2+P(x))\left(y^k+\sum_{j=0}^{k-1}\xi_{z,j}y^j\right)(1+O(x,y,\overline\xi))+t^\mu(F\big|_z+O(t,x,y,\overline\xi))\
.$$ Again the smoothness and transversality statement in Lemma
\ref{llz}(2) convert the condition to fit the deformation pattern
$h^{(2)}_z(x,y)$ and to realize a deformation of the union of
$kL(z)$ with the germs of $E$ at the points $E\cap L(z)$ into in
immersed cylinder, into a system of equations
\begin{equation}\xi_{z,j}+O(\overline\xi^2)=t^{\mu(k-j)/k}(b_{z,j}+O(t,\overline\xi)),\
1\le j<k,\quad\xi_{z,0}=t^{\mu}(-F\big|_z+O(t,\overline\xi))\
.\label{elz}\end{equation}

\smallskip

{\bf (3c)} Let $\Psi'(x,y)\in{\mathcal P}(L')$ be a deformation
pattern for the component $(L')^{s'}$ of $C$. In agreement with
(\ref{e40}) we can write (slightly modifying notations)
$$\Psi'(x,y)=b'+(y+x^2)f'(y),\quad f'(y)=y^{s'}+b'_{s'-1}y^{s'-1}+...+b'_0\
,$$ with $f'(y)$ defined as in Lemma \ref{l3}. On the other hand, in
the coordinates $x,y$ introduced in Step 1 of part (5), expression
(\ref{enn14}) with $\xi_0=t^\mu(c+O(t))$ reads as
$$(y+xy+x^2)\left(y^{s'}+\sum_{k=0}^{s'-1}\xi'_ky^k+O(\overline\xi^2)\right)(1+O(x,y,\overline\xi))+t^\mu(b'+O(t,\overline\xi))\ .$$
Again in view of the smoothness and transversality statement in
Lemma \ref{l3}(3), the conditions to fit the given deformation
pattern $\Psi'$ together with the demand to realize a deformation in
a neighborhood of the point $L'\cap E$ as described in Lemma
\ref{l1}(2vi), turn into a system of equations
\begin{equation}\xi'_k+O(\overline\xi^2)=t^{\mu(s'+1-k)/(s'+1)}(b'_k+O(t,\overline\xi)),\quad k=0,...,s'-1\
.\label{enn16}\end{equation} Similarly, for the component
$(L'')^{s''}$ of $C$ with a given deformation pattern $\Psi''$, we
obtain equations
\begin{equation}\xi''_k+O(\overline\xi^2)=t^{\mu(s''+1-k)/(s''+1)}(b''_k+O(t,\overline\xi)),\quad k=0,...,s''-1\
.\label{enn17}\end{equation}

\smallskip

{\bf (3d)} In the right-hand side of
$N-1\overset{(\ref{eparam})}{=}\sum_{C^{(i)}\ne
kL(z)}(I\gam^{(i)}-|\gam^{(i)}|)+\sum_{C^{(i)}=kL(z)}k+s'+s''$
equations (\ref{enn15}), (\ref{elj}), (\ref{enn16}), and
(\ref{enn17}), we have $N-1$ unknown functions. Observe that the
linearization of this unified system has a nondegenerate diagonal
form. Hence it is uniquely soluble as far as we take
$\xi_0=t^\mu(a+O(t))$. This provides a uniquely defined
parameterized branch $V$ of $V(\overline\bp,C)$ matching all the
given deformation patterns.

\smallskip

Had we taken $\xi_0=t^{r\mu}(c+O(t))$ with some $r>1$ in the
presented construction, we were coming (due to the uniqueness claim)
to the branch of $V(\overline\bp,C)$, geometrically coinciding with
$V$, with the parametrization obtained from that of $V$ via the
parameter change $t\mapsto t^r$.

The proof of Lemma \ref{l4} and Proposition \ref{p5}(2) is
completed. \proofend

\subsection{Recursive formula} \label{sec:formula}

The set $\Pic(\Sig,E) \times \Z \times \Z_+^\infty \times
\Z_+^\infty$ is a semigroup with respect to the operation
$$(D',g',\alp',\bet')+(D'',g'',\alp'',\bet'')=(D'+D'',g'+g''-1,\alp'+\alp'',\bet'+\bet'')\ .$$
Denote by $A(\Sig,E)\subset\Pic(\Sig,E) \times \Z \times \Z_+^\infty
\times \Z_+^\infty$ the subsemigroup generated by
\begin{itemize}\item all the quadruples $(D,g,\alp,\bet)$ with reduced,
irreducible $D$, $g\ge0$, and $\alp,\bet$ satisfying (\ref{en1}),
and \item the quadruples $(-s(K_\Sig+E),0,e_s,e_s)$ and
$(-s(K_\Sig+E),0,0,2e_s)$ for all $s\ge2$.
\end{itemize} Notice that one may have $g<0$ for elements
$(D,g,\alp,\bet)\in A(\Sig,E)$; in such a case we always assume
$N_\Sig(D,g,\alp,\bet)=0$.

\begin{theorem}\label{t1}
Given an element $(D,g,\alp,\bet)\in A(\Sig,E)$ such that
$$R_\Sig(D,g,\bet)>0\quad\text{and}\quad(D,g,\alp,\bet)\ne(-s(K_\Sig+E),0,0,2e_s),\ s>1\ ,$$ we have
$$N_\Sig(D,g,\alp,\bet) = \sum_{j \geq 1,\ \bet_j > 0} j \cdot
N_\Sig(D,g,\alp+e_j,\bet-e_j)$$
\begin{equation}
+ \sum \binom{\alp}{\alp^{(1)},\ldots,\alp^{(m)}}
\frac{(n-1)!}{n_1!...n_m!}\binom{k+3}{3} \prod_{i\in S}\left(
\binom{\bet^{(i)}}{\gam^{(i)}} I^{\gam^{(i)}}
N_\Sig(D^{(i)},g^{(i)},\alp^{(i)},\bet^{(i)}) \right),
\label{eq:main}\end{equation} where
$$n =R_\Sig(D,g,\bet),\quad n_i = R_\Sig(D^{(i)},g^{(i)},\bet^{(i)}),\quad   i =1,...,m\
,$$ $$S=\{i\in[1,m]\ :\
(D^{(i)},g^{(i)},\alp^{(i)},\bet^{(i)})\ne(-s(K_\Sig+E),0,e_s,e_s),\
s\ge1\}\ ,$$ and the second sum is taken
\begin{itemize}
\item over all integer $k\ge0$ and all splittings in $A(\Sig,E)$
\begin{equation} \label{eq:splittingsq}
(D-E+k(K_\Sig+E),g',\alp',\bet') = \sum_{i=1}^m
(D^{(i)},g^{(i)},\alp^{(i)},\bet^{(i)})\ ,\end{equation} such that
$$g^{(i)}\ge0\quad\text{and}\quad(D^{(i)},g^{(i)},\alp^{(i)},\bet^{(i)})\ne(-s(K_\Sig+E),0,0,se_2),\
s\ge 1,\ i=1,...,m\ ,$$ of all possible quadruples
$(D-E+k(K_\Sig+E),g',\alp',\bet') \in A(\Sig,E)$ with $k\ge0$,
satisfying
\begin{enumerate}
\item[(a)] $\alp' \leq \alp,\quad\bet' \geq \bet,\quad g - g' = \|\bet' - \bet\| + 1$,
\item[(b)] each summand $(D^{(i)},g^{(i)},\alp^{(i)},\bet^{(i)})$ with $n_i = 0$ and $\alp^{(i)}=0$ appears in
(\ref{eq:splittingsq}) at most once,
\end{enumerate}
\item over all splittings
\begin{equation} \label{eq:splittingsbeta}
\bet' = \bet + \sum_{i=1}^m \gam^{(i)},\quad\|\gam^{(i)}\|
> 0, \ 1 \leq i \leq m
\end{equation}
satisfying the restrictions $\bet^{(i)} \geq \gam^{(i)}$ for all $1
\leq i \leq m$,
\end{itemize}
and, finally, the second sum in (\ref{eq:main}) is factorized by
simultaneous permutations in both splittings (\ref{eq:splittingsq})
and (\ref{eq:splittingsbeta}).
\end{theorem}

{\bf Proof.} Immediately follows from from Propositions \ref{pdeg1},
and \ref{p5}. So, the first summand corresponds to the degeneration
described in Proposition \ref{pdeg1}(1), and its multiplicity $j$
was computed in Proposition \ref{p5}(1). The second sum corresponds
to the degeneration described in Proposition \ref{pdeg1}(2), and its
multiplicity was computed in Proposition \ref{p5}(2). Notice only
that the subtraction of $-k(K_\Sig+E)$ in (\ref{eq:splittingsq})
corresponds to all possible combinations $s'L'\cup s''L''$,
$s'+s''=k$, which, by Proposition \ref{p5}(2), altogether contribute
to the multiplicity $(k+1)+k\cdot 2+(k-1)\cdot
3+...+(k+1)=\binom{k+3}{k}$. \proofend

\smallskip

Note that the second sum in the right hand side of (\ref{eq:main})
becomes empty if $D - E$ is not effective.

\begin{theorem} \label{prop:rigidcases}
All the numbers $N_\Sig(D,g,\alp,\bet)$ with $(D,g,\alp,\bet) \in
A(\Sig,E)$ and $R_\Sig(D,g,\bet)>0$ are determined recursively by
formula (\ref{eq:main}) and the initial values
$N_\Sig(D,g,\alp,\bet)$ indicated in Proposition \ref{pini}.
\end{theorem}

{\bf Proof.} Straightforward. \proofend

In the important case of genus zero, all components of degenerate
curves in Proposition \ref{pdeg1} are rational, and the general
formula (\ref{eq:main}) reduces to the count of only genus zero
terms. Furthermore, it can be simplified so that the terms
corresponding to multiple covers will not explicitly appear.

Let us write for short $N_\Sig(D,\alp,\bet):=N_\Sig(D,0,\alp,\bet)$
and introduce the semigroup
$A_0(\Sig,E)\subset\Pic(\Sig,E)\times(\Z_+^{\infty})^2$, generated
by the triples $(D,\alp,\bet)$ with reduced, irreducible $D$ and
$I(\alp+\bet)=DE$.

\begin{corollary}\label{czero}
Given an element $(D,\alp,\bet)\in A_0(\Sig,E)$ such that
$R_\Sig(D,0,\bet)>0$, we have $$ N_\Sig(D,\alp,\bet) = \sum_{j \geq
1,\ \bet_j > 0} j \cdot N_\Sig(D,\alp+e_j,\bet-e_j)$$ $$ + \sum
\frac{2^{\|\bet^{(0)}\|}I^{\bet^{(0)}}}{\bet^{(0)}!}\binom{k+3}{3}\binom{\alp}{\alp^{(0)},\ldots,\alp^{(m)}}
\binom{n-1}{n_1,\ldots,n_m}$$ \begin{equation}\times \prod_{i=1}^m
\left( \binom{\bet^{(i)}}{\gam^{(i)}} I^{\gam^{(i)}}
N_\Sig(D^{(i)},\alp^{(i)},\bet^{(i)}) \right)\
,\label{eqrat}\end{equation} where
$$n =R_\Sig(D,0,\bet),\quad n_i = R_\Sig(D^{(i)},0,\bet^{(i)}),\quad   i =1,...,m, $$
and the second sum is taken
\begin{itemize} \item over
all integers $k\ge0$ and vectors
$\alp^{(0)},\bet^{(0)}\in\Z_+^\infty$ such that $\alp^{(0)}\le\alp$
and $\bet^{(0)}\le\beta$;
\item over all sequences
\begin{equation}(D^{(i)},\alp^{(i)},\bet^{(i)})\in
A_0(\Sig,E),\quad 1\le i\le m \
,\label{eq:splittings0}\end{equation} such that
\begin{enumerate}
\item[(i)] $D^{(i)}\ne -(K_\Sig+E)$ and
$R(D^{(i)},0,\bet^{(i)})\ge0$ for all $i=1,...,m$,
\item[(ii)]
$D-E=\sum_{i=1}^mD^{(i)}-(k+I\alp^{(0)}+I\bet^{(0)})(K_\Sig+E)$,
\item[(iii)] $\sum_{i=0}^m\alp^{(i)}\le\alp,\quad\sum_{i=0}^m\bet^{(i)}\ge\bet$,
\item[(iv)] each
triple $(D^{(i)},0,\bet^{(i)})$ with $n_i=0$ appears in
(\ref{eq:splittings0}) at most once,
\end{enumerate}
\item over all sequences \begin{equation}\gam^{(i)}\in\Z_+^\infty,\quad
\|\gam^{(i)}\|=1,\quad i=1,...,m\ ,\label{eq:splittings1}
\end{equation} satisfying
$$\bet^{(i)}\ge\gam^{(i)}, \ i=1,...,m,\quad\text{and}\quad
\sum_{i=1}^m\left(\bet^{(i)}- \gam^{(i)}\right)=\bet-\bet^{(0)}\ ,$$
\end{itemize}
and the second sum in (\ref{eqrat}) is factorized by simultaneous
permutations in the sequences (\ref{eq:splittings0}) and
(\ref{eq:splittings1}).
\end{corollary}

\section{Counting curves on $\PP^2_7$}

Since $\PP^2_{6,1}$ is almost Fano in the sense of \cite[Section
4.1]{Va}, i.e. $K_{\PP^2_{6,1}}$ is negative on all curves except
for $E$. Hence \cite[Theorem 4.2]{Va} yields the following relation
between $GW_g(\PP^2_7,D)$ and enumerative invariants of
$\PP^2_{6,1}$ (abusing notations, we use the same symbol for divisor
classes corresponding via a natural isomorphism
$\Pic(\PP^2_7)\simeq\Pic(\PP^2_{6,1})$).

\begin{theorem}\label{gw}
For any effective divisor class $D$ on $\PP^2_7$, the Gromov-Witten
invariant $\qquad$ $GW_g(\PP^2_7,D)$ is given by the formula
\begin{equation}GW_g(\PP^2_7,D)=\sum_{i\ge0}\binom{DE+2i}{i}N_{\PP^2_{6,1}}(D-iE,g,0,(DE+2i)e_1)\
. \label{eGW}
\end{equation}
\end{theorem}

{\bf Proof}. We have to make only one remark. The maps $\bn:C'\cup
C''\to\PP^2_{6,1}$. which are counted in \cite[Theorem 2]{Va}, where
the components of $C'$ are isomorphically mapped onto $E$, and the
components of $C''$ are not mapped onto $E$, the part $C''$ is
connected. Thus, by Proposition \ref{pn1}(2) it is smooth of genus
$g$, and its image in $\PP^2_{6,1}$ crosses $E$ transversally. The
binomial coefficient reflects the choice of intersection points of
$C'$ and $C''$ among the points of $\bn^{-1}(\bn_*(C'')\cap E)$.
\proofend

\begin{example}
For the popular case of $D=-2K_{\PP^2_6}=6L-2(E_1+...+E_6)$ (cf.
\cite[Page 119, Table II]{DFI}, \cite[Page 88]{GP}, \cite[Section
9.2]{Va}), we compute Gromov-Witten invariants for all genera in two
ways, once using Theorem \ref{t1} for $\PP^2_{5,1}=\PP^2_6$, and
another time using Theorem \ref{t1} for $\PP^2_{6,1}$ with
$D=6L-2(E_2-...-E_7)$ and then applying Theorem \ref{gw} (cf.
\cite[Section 9.2]{Va}):
\begin{center}
\begin{tabular}{|l||c|c|c|c|c|c|}
\hline $g$ &\makebox[0.35cm] 0 &\makebox[0.35cm] 1 &\makebox[0.35cm]
2&\makebox[0.35cm] 3 &\makebox[0.35cm] 4
  \\
\hline
$GW$ & 3240 & 1740 & 369 & 33 & 1  \\
\hline
\end{tabular}
\end{center}
For the case of $D=-2K_{\PP^2_7}=6L-2(E_1+...+E_7)$ (cf. \cite[Page
88]{GP}), our formula gives the following values for the
Gromov-Witten invariants:
\begin{center}
\begin{tabular}{|l||c|c|c|c|c|}
\hline $g$ &\makebox[0.35cm] 0 &\makebox[0.35cm] 1 &\makebox[0.35cm]
2&\makebox[0.35cm] 3
  \\
\hline
$GW$ & 576 & 204 & 26 & 1  \\
\hline
\end{tabular}
\end{center}
\end{example}

\vskip10pt

{\ncsc Institute of Mathematics, \\[-21pt]

Hebrew University of Jerusalem, \\[-21pt]

Givat Ram, Jerusalem, 91904, Israel} \\[-21pt]

\emph{E-mail}: {\ntt mendy.shoval@mail.huji.ac.il}

\vskip10pt

{\ncsc School of Mathematical Sciences, \\[-21pt]

Raymond and Beverly Sackler Faculty of Exact Sciences,\\[-21pt]

Tel Aviv University,
Ramat Aviv, 69978 Tel Aviv, Israel} \\[-21pt]

\emph{E-mail}: {\ntt shustin@post.tau.ac.il}

\end {document}